\documentclass[12pt,leqno]{amsart}
\usepackage{amssymb}
\usepackage{amsmath,amssymb}
\newtheorem{theorem}{Theorem}[section]
\newtheorem{corollary}{Corollary}[section]

\newtheorem{proposition}{Proposition}[section]

\newcommand{\dddd}{\partial_z}
\newcommand{\ddddd}{\partial_{\overline{z}}}
\begin{document}
\theoremstyle{plain}
\newtheorem{MainThm}{Theorem}
\newtheorem{thm}{Theorem}[section]
\newtheorem{clry}[thm]{Corollary}
\newtheorem{prop}[thm]{Proposition}
\newtheorem{lem}[thm]{Lemma}
\newtheorem{deft}[thm]{Definition}
\newtheorem{hyp}{Assumption}
\newtheorem*{ThmLeU}{Theorem (J.~Lee, G.~Uhlmann)}

\theoremstyle{definition}
\newtheorem{rem}[thm]{Remark}
\newtheorem*{acknow}{Acknowledgments}
\numberwithin{equation}{section}
\newcommand{\eps}{{\varphi}repsilon}
\renewcommand{\d}{\partial}
\newcommand{\re}{\mathop{\rm Re} }
\newcommand{\im}{\mathop{\rm Im}}
\newcommand{\R}{\mathbf{R}}
\newcommand{\C}{\mathbf{C}}
\newcommand{\N}{\mathbf{N}}
\newcommand{\D}{C^{\infty}_0}
\renewcommand{\O}{\mathcal{O}}
\newcommand{\dbar}{\overline{\d}}
\newcommand{\supp}{\mathop{\rm supp}}
\newcommand{\abs}[1]{\lvert #1 \rvert}
\newcommand{\csubset}{\Subset}
\newcommand{\detg}{\lvert g \rvert}
\title[partial Dirichlet-to-Neumann map]
{Global Uniqueness from partial Cauchy data
in two dimensions}

\author[O. Imanuvilov]{Oleg Yu. Imanuvilov}
\address{Department of Mathematics, Colorado State
University, 101 Weber Building, Fort Collins CO, 80523 USA\\
e-mail: oleg@math.colostate.edu}
\thanks{First author partly supported by NSF grant DMS 0808130}

\author[G. Uhlmann]{Gunther Uhlmann}
\address{Department of Mathematics, University of Washington, Seattle,
WA 98195 USA\\
e-mail: gunther@math.washington.edu}
\thanks{Second author partly supported by NSF and a Walker Family Endowed
Professorship}

\author[M. Yamamoto]{Masahiro Yamamoto}
\address{Department of Mathematics, University of Tokyo, Komaba, Meguro,
Tokyo 153, Japan \\e-mail:  myama@ms.u-tokyo.ac.jp}

\begin{abstract}
We prove for a two dimensional bounded  domain
that the Cauchy data for the Schr\"odinger equation measured on an
arbitrary open subset of the boundary determines uniquely the potential.
This implies, for the conductivity equation, that if we measure the
current fluxes at the boundary on an arbitrary open subset of the
boundary produced by voltage potentials supported in the same subset,
we can determine uniquely the conductivity.
We use Carleman estimates with degenerate
weight functions to construct appropriate complex geometrical
optics solutions to prove the results.
\end{abstract}

\maketitle \setcounter{tocdepth}{1} \setcounter{secnumdepth}{2}

\section{\bf Introduction}

We consider the problem of determining a complex-valued potential $q$ in
a bounded  two dimensional domain from the Cauchy data
measured on an arbitrary open subset of the boundary for the associated
Schr\"odinger equation $\Delta + q$.
A motivation comes from the classical inverse problem of
electrical impedance tomography
problem. In this inverse problem one attempts to
determine the electrical conductivity of a body by measurements of
voltage and current on the boundary of the body. This problem was
proposed by Calder{\'o}n \cite{C} and is also known as Calder\'on's
problem.
In dimensions $n\ge 3$, the first global uniqueness result for
$C^2$-conductivities was proven in \cite{SU}. In
\cite{PPU}, \cite{BT} the global uniqueness result was extended to less regular
conductivities. Also see \cite{GLU} as for the determination of more
singular conormal conductivities. In dimension $n\ge 3$ global
uniqueness was shown for the Schr\"odinger equation with bounded
potentials in \cite{SU}.
The case of more singular conormal potentials was studied
in \cite{GLU}.

In two dimensions the first global uniqueness result for
Calder\'on's problem was obtained in \cite{N} for
$C^2$-conductivities. Later the regularity assumptions were relaxed
in  \cite{B-U}, and \cite{AP}. In particular, the paper
\cite{AP} proves
uniqueness for $L^\infty$- conductivities. In two dimensions a recent
result of Bukgheim \cite{Bu} gives unique identifiability of the
potential from Cauchy data measured on the whole boundary for the associated Schr\"odinger
equation.
As for the uniqueness in determining two coefficients, see
\cite{ChengYama}, \cite{KU}.

In all the above mentioned articles, the measurements are made on
the whole boundary. The purpose of this paper is to show the global
uniqueness in two dimensions, both for the Schr\"odinger and
conductivity equation, by measuring all the Neumann data on an
arbitrary open subset $\widetilde{\Gamma}$ of the boundary produced
by inputs of Dirichlet data supported on $\widetilde{\Gamma}$. We
formulate this inverse problem more precisely below. Let
$\Omega\subset \R^2$ be a  bounded domain with smooth boundary, and
let $\nu$ be the unit outward normal vector to $\partial\Omega$. We
denote $\frac{\partial u}{\partial\nu} = \nabla u \cdot\nu$.  A
bounded and non-zero function $\gamma(x)$ (possibly complex-valued)
models the electrical conductivity of $\Omega$. Then a potential
$u\in H^1(\Omega)$ satisfies the Dirichlet problem
\begin{equation}\label{eq:0.2}
\begin{array}{rcl}
\mbox{div}(\gamma\nabla u) & = & 0 \mbox{ in }\Omega, \\
u\big|_{\partial\Omega} & = & f,
\end{array}
\end{equation}
where $f\in H^{\frac{1}{2}}(\partial\Omega)$ is a given boundary
voltage potential. The Dirichlet-to-Neumann (DN) map is defined by
\begin{equation}\label{eq:0.3}
\Lambda_\gamma(f)= \gamma\frac{\partial u}{\partial \nu
}\Big|_{\partial\Omega}.
\end{equation}

This problem can be reduced to studying the set of Cauchy data for
the Schr\"odinger equation  with the potential $q$ given by:
\begin{equation}\label{eq:2.2}
q=\frac{\Delta\sqrt{\gamma}}{\sqrt{\gamma}}.
\end{equation}

More generally we define the set of Cauchy data for a bounded potential $q$ by:
\begin{equation}
\widehat{C_q} =
\left\{\left(u|_{\partial\Omega}, \frac{\partial u}{\partial
\nu}\Big|_{\partial\Omega}\right)\mid (\Delta+q) u= 0\hbox{ on
}\Omega,\,\, \ u\in H^1(\Omega)\right\}.
\end{equation}
We have $\widehat{C_q}\subset H^{\frac{1}{2}}(\partial\Omega)\times
H^{-\frac{1}{2}}(\partial\Omega)$.

Let $\widetilde \Gamma \subset
\partial\Omega$ be a non-empty open subset of the boundary.  Denote
$\Gamma_0=\partial\Omega\setminus \widetilde \Gamma.$

Our main result gives global
uniqueness by measuring the Cauchy data on $\widetilde \Gamma$.
Let $q_j \in C^{1+\alpha}(\overline \Omega)$, $j=1,2$ for some
$\alpha>0$ and let $q_j$ be complex-valued.
Consider the following sets of Cauchy data on an $\widetilde \Gamma$:
\begin{equation}
\mathcal{C}_{q_j}=\left\{\left(u|_{\widetilde\Gamma}, \frac{\partial
u}{\partial \nu}\Big|_{\widetilde\Gamma}\right)\mid (\Delta+q_j) u=
0\hbox{ on }\Omega,\, u\vert_{\Gamma_0}=0,\,\, u\in
H^1(\Omega)\right\}, \quad j=1,2.
\end{equation}
\begin{theorem}\label{main}
Assume $ \mathcal{C}_{q_1}= \mathcal{C}_{q_2}. $ Then $
q_1\equiv q_2. $
\end{theorem}

Using Theorem~\ref{main} one concludes immediately as
a corollary the following global identifiability result
for the conductivity
equation (\ref{eq:0.2}).

\begin{corollary}\label{coro}
With some $\alpha > 0$, let $\gamma_j \in
C^{3+\alpha}(\overline\Omega)$, $j=1,2$, be non-vanishing functions.
Assume that
$$
{\Lambda}_{\gamma_1}(f)={\Lambda}_{\gamma_2}(f)\mbox{ in }
\widetilde\Gamma \mbox{ for all } f\in H^{1\over 2}(\Gamma
),\,\,\mbox{supp}\thinspace f\subset\widetilde\Gamma.
$$
Then $\gamma_1=\gamma_2$.
\end{corollary}
Notice that Theorem \ref{main}  does not assume that $\Omega$ is simply connected. An interesting inverse problem is  where one can determine the potential or conductivity in a region of the plane with holes by measuring  the Cauchy data only on the accessible boundary. Let $\Omega,D$ be domains in $\R^2$ with smooth boundary such that $\overline D\subset \Omega.$ Let $V\subset \partial\Omega$ be an open set.  Let $q_{j}\in C^{1+\alpha}(\overline{\Omega\setminus D}),$  for some $\alpha>0$ and  $j=1,2.$ Let us consider the following set of partial Cauchy data
$$
\tilde C_{q_j}=\{(u\vert_V,\frac{\partial u}{\partial \nu}\vert_V)\vert (\Delta+q_j)u=0\,\,\mbox{in}\,\,\Omega\setminus\overline D, \, u\vert_{\partial D\cup \partial\Omega\setminus V}=0, u\in H^1(\Omega)\}.
$$
\begin{corollary}\label{coro1} Assume $\tilde C_{q_1}=\tilde C_{q_2}$. Then $q_1=q_2.$
\end{corollary}

A similar result holds for the conductivity equation.
\begin{corollary}\label{coro2} Let $\gamma_j\in C^{3+\alpha}(\overline{\Omega\setminus D})$ $j=1,2$ be non vanishing functions. Assume
$$
\Lambda_{\gamma_1}(f)=\Lambda_{\gamma_2}(f)\,\, \mbox{in}\,\,\, V\quad \forall f\in H^\frac 12(\partial(\overline\Omega\setminus D)),\,\,\mbox{supp}\, f\subset V
$$ Then  $\gamma_1=\gamma_2$.
\end{corollary}

Another application of Theorem\ref{main} is to the anisotropic conductivity problem. In this case the conductivity depends on direction  and is represented by a positive definite symmetric matrix
$$
\sigma=\{\sigma^{ij}\}\quad\mbox{on}\,\,\Omega.
$$
The conductivity equation with voltage potential $g$ on $\partial\Omega$ is given by
$$
\sum_{i,j=1}^2\frac{\partial}{\partial x_i}(\sigma^{ij}\frac{\partial u}{\partial x_j})=0\quad\mbox{in}\,\,\Omega,
$$
$$
u\vert_{\partial\Omega}=g.
$$
The Dirichlet-to-Neumann map is defined  by
$$
\Lambda_{\sigma}(g)=\sum_{i,j=1}^2\sigma^{ij}\nu_i\frac{\partial u}{\partial x_j}\vert_{\partial\Omega}.
$$
It has been known for a long time that $\Lambda_\sigma$ does not determine $\sigma$ uniquely in the anisotropic case \cite{K-V}.
Let $F:\overline\Omega\rightarrow \overline\Omega$ be a diffeomorphism such that $F(x)=x$ for and $x$ from $\partial\Omega.$ Then
$$
\Lambda_{F_*\sigma}=\Lambda_\sigma,
$$
where
\begin{equation}\label{star}F_*\sigma=\frac{(DF)\cdot\sigma\cdot(DF)^T\cdot F^{-1}}{\vert det DF\vert}.
\end{equation}
Here $DF$ denotes the differential of $F,$ $(DF)^T$ its transpose and the composition in (\ref{star}) is matrix composition.  The question of whether one can determine the conductivity up to the obstruction (\ref{star}) has been solved  for $C^2$ conductivities in \cite{N},  Lipschitz conductivities in \cite{SuU} and merely $L^\infty$ conductivities in \cite{ALP}.  The method of proof in all these papers is the reduction to the isotropic case performed using isothermal coordinates \cite{S}. Using the same method and Corollary \ref{coro} we obtain the following result
\begin{theorem}Let $\sigma_k=\{\sigma_k^{ij}\}\in C^{3+\alpha}(\overline\Omega)$ for $k=1,2$ and some positive $\alpha.$ Suppose that $\sigma_k$ are positive definite symmetric matrices on $\overline\Omega$. Let $\Gamma\subset \partial\Omega$ be some open set. Assume
$$
\Lambda_{\sigma_1}(g)\vert_\Gamma=\Lambda_{\sigma_2}(g) \vert_\Gamma\quad \forall g\in H^\frac 12(\partial\Omega),\,\mbox{supp}\, g\subset \Gamma.
$$ Then there exists a diffeomorphism
$$
F:\overline\Omega\rightarrow \overline\Omega, \quad F\vert_{\partial\Omega}=\mbox{Identity},\quad F\in C^{4+\alpha}(\bar\Omega), \alpha >0$$
such that
$$
F_*\sigma_1=\sigma_2.
$$
\end{theorem}

To the authors' knowledge, there are no uniqueness results similar to
Theorem \ref{main} with Dirichlet data supported and Neumann  data measured
on the same arbitrary open subset of the boundary, even for smooth
potentials or conductivities. In dimension $n\ge 3$ Isakov \cite {I} proved global uniqueness assuming that $\Gamma_0$ is a subset of a plane or a sphere.
In dimensions $n\ge 3$, \cite{BuU} proves global uniqueness in
determining a bounded potential for the Schr\"odinger equation
with Dirichlet data supported on the whole boundary and
Neumann data measured in roughly half the boundary.
The proof relies on a Carleman estimate with a linear
weight function. This implies a similar result for the conductivity
equation with $C^2$ conductivities. In \cite{K} the regularity
assumption on the
conductivity was relaxed to $C^{3/2+\alpha}$ with some $\alpha>0$. The
corresponding stability estimates are proved in \cite{HW}. As for
the stability estimates for the magnetic Schr\"odinger equation with
partial data, see \cite{T}.  In \cite{KSU}, the result in \cite{BuU}
was generalized to show that by all possible pairs of Dirichlet data
on an arbitrary open subset $\Gamma_+$ of the boundary and Neumann
data on a slightly larger open domain than $\partial\Omega\setminus
\Gamma_+$, one can uniquely determine the potential.  The method
of the proof uses Carleman estimates with non-linear weights. The case
of the magnetic Schr\"odinger equation was considered in \cite{DKSjU}
and an improvement on the regularity of the coefficients is done in
\cite{KS}.
Stability estimates for the magnetic Schr\"odinger equation with
partial data were proven in \cite{T}.

In two dimensions the first general result was given by the authors
in \cite{IGM}. It is shown that the same global
uniqueness result as \cite{KSU} holds in this case.
The two dimensional case has special features since one can construct
a much larger set of complex geometrical optics solutions than in higher
dimensions.  On the other hand, the problem is formally
determined in two dimensions and therefore more difficult.
The proof of our main result
\cite{IGM} is based on the construction of appropriate complex
geometrical optics solutions by Carleman
estimates with degenerate weight functions.

This paper is composed of four sections.  In Section 2, we establish
our key Carleman estimates, and in Section 3, we construct complex
geometrical optics solutions.  In Section 4, we complete the proof
of Theorem \ref{main}. In the Appendix we
consider the solvability of the Cauchy Riemann equations with Cauchy data
 on a subset of the boundary. We also establish a Carleman estimate for
 Laplace's equation with degenerate harmonic weights.

\section{\bf Carleman estimates with degenerate weights}

Throughout the paper we use the following notations:
\\

\noindent {\bf Notations}

$i=\sqrt{-1}$, $x_1, x_2, \xi_1, \xi_2 \in
\R$, $z=x_1+ix_2$, $\zeta=\xi_1+i\xi_2$, $\overline{z}$ denotes
the complex conjugate of $z \in \C$.  We identify $x = (x_1,x_2)
\in \R^2$ with $z = x_1 +ix_2 \in \C$.
$\dddd = \frac 12(\partial_{x_1}-i\partial_{x_2})$,
$\ddddd = \frac12(\partial_{x_1}+i\partial_{x_2}),$
$H^{1,\tau}(\Omega)$ denotes
the space $H^1(\Omega)$ with norm $\Vert
v\Vert^2_{H^{1,\tau}(\Omega)}=\Vert
v\Vert^2_{H^{1}(\Omega)}+\tau^2\Vert v\Vert^2_{L^2(\Omega)}.$ The
tangential derivative on the boundary is given by
$\partial_{\vec\tau}=\nu_2\frac{\partial}{\partial x_1}
-\nu_1\frac{\partial}{\partial x_2},$ with $\nu=(\nu_1, \nu_2)$ the
unit outer normal to $\partial\Omega,$ $B(\widehat x,\delta)=\{x\in
\R^2\vert \vert x-\widehat x\vert< \delta\},$
$f(x):\R^2\rightarrow \R^1$, $f''$ is the Hessian matrix  with entries
$\frac{\partial^2 f}{\partial x_i\partial x_j}.$ $\mathcal L(X,Y)$
denotes the Banach space of all bounded linear operators from a Banach
space $X$ to another Banach space $Y$.

Let $\Phi(z)=\varphi(x_1,x_2)+i\psi(x_1,x_2) \in C^2(\overline{\Omega})$
be a holomorphic function in $\Omega$ with real-valued $\varphi$ and
$\psi$:
\begin{equation}\label{zzz}
\ddddd\Phi(z) = 0\quad \mbox{in} \,\,\Omega.
\end{equation}
Denote by $\mathcal H$ the set of critical points of a function
$\Phi$
$$
\mathcal H = \{z\in\overline\Omega\vert \dddd \Phi(z)=0 \}.
$$
Assume that $\Phi$ has no critical points on the boundary, and
that all the critical points are nondegenerate:
\begin{equation}\label{mika}
\mathcal H\cap \partial\Omega=\{\emptyset\},\quad
\partial_z^2\Phi(z)\ne 0, \quad \forall z\in \mathcal H.
\end{equation}
Then we know that $\Phi$  has only a finite number of critical points
and we can set:
\begin{equation}\label{mona}
\mathcal H = \{ \widetilde{x}_1, ..., \widetilde{x}_{\ell} \}.
\end{equation}

Consider the following problem
\begin{equation}\label{(2.26)}
\Delta u+q_0u=f\quad \mbox{in}\,\,\Omega,\quad u\vert_{\Gamma_0}=g,
\end{equation}
where $\nu$ is the unit outward normal vector to $\partial\Omega$ and
$$
\Gamma_0=\{x\in \partial\Omega\vert (\nu,\nabla\varphi)=0\}.
$$
We have
\begin{proposition}\label{Proposition 2.3}
Let $q_0\in L^\infty(\Omega).$ There exists $\tau_0>0$ such that for
all $\vert\tau\vert>\tau_0$  there exists a solution to problem
(\ref{(2.26)}) such that \begin{equation} \label{(2.27)}\Vert
ue^{-\tau\varphi}\Vert_{L^2(\Omega)}\le C(\Vert
fe^{-\tau\varphi}\Vert_{L^2(\Omega)}/\root\of{\vert \tau\vert}+\Vert
ge^{-\tau\varphi}\Vert_{L^2(\Gamma_0)}).
\end{equation}
\end{proposition}

For the proof, see Proposition 2.2 in \cite{IGM} and Proposition
\ref{Theorem 2.1} in appendix.

Let us introduce the operators:
$$
\partial_{\overline z}^{-1}g=\frac{1}{2\pi i}\int_\Omega
\frac{g(\zeta,\overline \zeta)}{\zeta-z} d\zeta\wedge d\overline\zeta=-\frac
1\pi\int_\Omega \frac{g(\zeta,\overline\zeta)}{\zeta-z}d\xi_2d\xi_1,
$$
$$
\partial_{ z}^{-1}g=-\frac{1}{2\pi i}\overline{\int_\Omega
\frac{\overline g(\zeta,\overline\zeta)}{\zeta-z} d\zeta\wedge
d\overline\zeta}=-\frac 1\pi\int_\Omega
\frac{g(\zeta,\overline\zeta)}{\overline\zeta-\overline z}d\xi_2d\xi_1
= \overline{\partial^{-1}_{\overline{z}}\overline{g}}.
$$
See e.g., pp.28-31 in \cite{VE} where $\ddddd^{-1}$ and $\dddd^{-1}$
are denoted by $T$ and $\overline{T}$ respectively.
Then we know (e.g., p.47 and p.56 in \cite{VE}):
\begin{proposition}\label{Proposition 3.0}
{\bf A)} Let $m\ge 0$ be an integer number and $\alpha\in (0,1).$ The
operators $\partial_{\overline z}^{-1},\partial_{ z}^{-1}\in
\mathcal L(C^{m+\alpha}(\overline{\Omega}),C^{m+\alpha+1}
(\overline{\Omega})).$
\newline
{\bf B}) Let $1\le p\le 2$ and $ 1<\gamma<\frac{2p}{2-p}.$ Then
 $\partial_{\overline z}^{-1},\partial_{ z}^{-1}\in
\mathcal L(L^p( \Omega),L^\gamma(\Omega)).$
\end{proposition}

We define two other operators:
\begin{equation}\label{(3.1)}
R_{\Phi,\tau}g=e^{\tau(\overline {\Phi(z)}
- {\Phi(z)})}\partial_{\overline z}^{-1}(ge^{\tau(
{\Phi(z)}-\overline {\Phi(z)})}),\quad \widetilde
R_{\Phi,\tau}g=e^{\tau(\overline {\Phi(z)}-{\Phi(z)})}\partial_{
z}^{-1}(ge^{\tau( {\Phi(z)}-\overline {\Phi(z)})}).
\end{equation}

\begin{proposition}\label{Proposition 3.1}
Let $g\in
C^\alpha(\overline\Omega)$ for some positive $\alpha.$ The function
$R_{\Phi,\tau}g$ is a solution to
\begin{equation}\label{(3.2)}
\partial_{\overline z}R_{\Phi,\tau}g
- \tau(\overline{\dddd \Phi(z)})R_{\Phi,\tau}g
= g \quad\mbox{in}\,\,\Omega.
\end{equation}
The function $\widetilde R_{\Phi,\tau}g$ solves
\begin{equation}\label{(3.3)}
\partial_{ z}\widetilde R_{\Phi,\tau}g
+ \tau (\dddd\Phi(z))\widetilde R_{\Phi,\tau}g
=g\quad\mbox{in}\,\,\Omega.
\end{equation}
\end{proposition}

The proof is done by direct computations (see the proof of
Proposition 3.3 in \cite{IGM}).

Denote
$$
\mathcal O_\epsilon=\{ x\in \Omega\vert dist(x,\partial\Omega)
\le \epsilon\}.
$$

\begin{proposition}\label{Proposition 3.2}
Let $g\in C^1(\Omega)$ and $g\vert_{\mathcal O_\epsilon}=0$,
$g(x)\ne 0 $ for all $x\in\mathcal H.$ Then
\begin{equation}\label{(3.4)}
\vert R_{\Phi,\tau}g(x)\vert
+ \vert\widetilde R_{\Phi,\tau} g(x)\vert\le C \max_{x\in\mathcal
H}\vert g(x) \vert/\tau
\end{equation}
for all $x\in \mathcal O_{\epsilon/2}.$
If $g\in C^2(\overline \Omega)$ and $g\vert_{\mathcal H}=0$, then
\begin{equation}\label{(3.5)}\vert R_{\Phi,\tau} g(x)\vert+\vert
\widetilde R _{\Phi,\tau}g(x)\vert\le C /\tau^2\end{equation} for
all $x\in \mathcal O_{\epsilon/2}.$
\end{proposition}

The proof uses the Cauchy-Riemann equations and
stationary phase (e.g., Section 4.5.3 in \cite{E}, Chapter VII, \S 7.7
in \cite{H}).
See also the proof of Proposition 3.4 in \cite{IGM}.

Denote
$$
r(z)=\Pi_{k=1}^\ell(z-\widetilde{z}_k)\,\,\,\mbox{where} \,\,\mathcal
H=\{\widetilde{x}_1,\dots, \widetilde{x}_\ell\}, \,\,\widetilde{z}_k=\widetilde{x}_{1,k}+i\widetilde{x}_{2,k}.
$$

The following proposition can be proved similarly to Proposition 3.5
in \cite{IGM}:
\begin{proposition}\label{Proposition 3.3}
Let $g\in C^1(\overline\Omega)$ and $g\vert_{\mathcal O_\epsilon}=
0.$ Then for each $\delta\in (0,1)$, there exists a constant
$C(\delta) > 0$ such that
\begin{equation}\label{(3.6)}
\Vert \widetilde R_{\Phi,\tau} (\overline
{r(z)}g)\Vert_{L^2(\Omega)} \le C(\delta)\Vert
g\Vert_{C^1(\overline\Omega)}/\vert\tau\vert^{1-\delta},\quad \Vert
R_{\Phi,\tau} (r(z)g)\Vert_{L^2(\Omega)}\le C(\delta)\Vert
g\Vert_{C^1(\overline\Omega)}/\vert\tau\vert^{1-\delta}.
\end{equation}
\end{proposition}
Henceforth we set $\psi_1 \equiv\mbox{Re} \partial_z\Phi
= \partial_{x_1}\varphi$ and
$\psi_2 \equiv \mbox{Im}\partial_z\Phi = \partial_{x_1}\psi$.
We also need the following result, which we can be proven
in the same way as  Proposition 2.1 in \cite{IGM}. Note that
$$
\partial_{x_1}(e^{-i\tau\psi}\widetilde{v})e^{i\tau\psi}
= \partial_{x_1}\widetilde{v} - i\tau\psi_2\widetilde{v}
$$
and
$$
\partial_{x_2}(e^{-i\tau\psi}\widetilde{v})e^{i\tau\psi}
= \partial_{x_2}\widetilde{v} - i\tau\psi_1\widetilde{v},
$$
etc. which follow from the Cauchy-Riemann equations.
\begin{proposition} \label{Proposition 2.1}
Let $\Phi$ satisfy (\ref{zzz}) and (\ref{mika}). Let $\widetilde
f\in L^2(\Omega)$ and $\widetilde v$ be solution to
\begin{equation}\label{zina}
2 \partial_z\widetilde v -\tau(\partial_z\Phi)\widetilde v
=\widetilde f\quad \mbox{in }\,\Omega
\end{equation}
or $\widetilde v$ be solution to
\begin{equation}\label{zina1}
2\partial_{\overline z}v - \tau(\partial_{\overline{z}}
\overline\Phi)\widetilde v
=\widetilde f\quad\mbox{ in }\,\Omega.
\end{equation}
In the case (\ref{zina}) we have
\begin{eqnarray}\label{vika1}
\Vert \partial_{x_1}(e^{-i\tau\psi}\widetilde{v})\Vert^2_{L^2(\Omega)}
- \tau\int_{\partial\Omega}(\nabla\varphi,\nu)\vert
\widetilde v\vert^2d\sigma\nonumber\\
+ \mbox{Re}\int_{\partial\Omega}i\left(\left(\nu_2
\frac{\partial}{\partial x_1}-\nu_1 \frac{\partial}{\partial
x_2}\right) \widetilde v\right)\overline{\widetilde v}d\sigma
+ \Vert \partial_{x_2}(e^{-i\tau\psi}\widetilde{v})\Vert^2
_{L^2(\Omega)}
= \Vert\widetilde f\Vert^2_{L^2(\Omega).}
\end{eqnarray}
In the case that $\widetilde v$ solves (\ref{zina1}) we have
\begin{eqnarray} \label{vika2}
\Vert \partial_{x_1}(e^{i\tau\psi}\widetilde{v})\Vert_{L^2(\Omega)}
- \tau\int_{\partial\Omega}(\nabla\varphi,\nu)\vert \widetilde
v\vert^2d\sigma + \mbox{Re}\int_{\partial\Omega}i\left(\left(-\nu_2
\frac{\partial}{\partial x_1}+\nu_1 \frac{\partial}{\partial
x_2}\right)\widetilde v\right)\overline{\widetilde
v}d\sigma\nonumber\\
+ \Vert \partial_{x_2}(e^{i\tau\psi}\widetilde{v})\Vert^2
_{L^2(\Omega)} = \Vert\widetilde f\Vert^2_{L^2(\Omega)}.
\end{eqnarray}
\end{proposition}

We have
\begin{proposition}\label{Proposition 3.22}
Let $g\in C^2(\Omega), g\vert_{\mathcal O_\epsilon} = 0$ and
$g\vert_{\mathcal H}= 0$.  Then
\begin{equation}\label{(3.4A)}
\left\Vert R_{\Phi,\tau}g+\frac {g}{\tau\overline{\partial_z \Phi}}
\right\Vert_{L^2(\Omega)}
+ \left\Vert\widetilde R_{\Phi,\tau} g-\frac {g}{\tau{\partial_z
\Phi}}\right\Vert_{L^2(\Omega)}
= o\left(\frac 1\tau\right)\quad \mbox{as}\,\,
\vert\tau\vert\rightarrow +\infty.
\end{equation}
\end{proposition}
\begin{proof}
By (\ref{mika}) and Proposition \ref{Proposition 3.2}
\begin{equation}\label{opa}
\Vert \widetilde{R}_{\Phi,\tau} g\Vert _{C(\overline{\mathcal
O_{\frac \epsilon 2}})}+ \Vert R_{\Phi,\tau}
g\Vert_{C(\overline{\mathcal O_{\frac \epsilon 2}})} = o\left(\frac
{1}{\tau}\right).
\end{equation}
Therefore instead of  (\ref{(3.4A)}) it suffices to prove
\begin{equation}\label{(3.4AA)}
\left\Vert \chi_1 R_{\Phi,\tau}g
+\frac {g}{\tau\overline{\partial_z \Phi}}\right\Vert_{L^2(\Omega)}
+ \left\Vert \chi_1 \widetilde R_{\Phi,\tau} g-\frac {g}{\tau{\partial_z
\Phi}}\right\Vert_{L^2(\Omega)}
=o\left(\frac 1\tau\right)\quad \mbox{as}\,\,
\vert\tau\vert\rightarrow +\infty,
\end{equation}
where $\chi_1\in C^\infty_0(\Omega)$ and
$\chi_1\vert_{\Omega\setminus\mathcal O_{\epsilon/2}}\equiv 1.$  Denote
$w=\chi_1\widetilde R_{\Phi,\tau}g-\frac {g}{\tau{\partial_z
\Phi}}.$  Here we note that $\frac{g}{\dddd\Phi} \in
L^\infty({\Omega})$.  This follows from (\ref{mika}),
$g \in C^1(\overline{\Omega})$ and $g\vert_{\mathcal H} = 0.$
Then (2.8) and $g\vert_{{\mathcal O}_{\varepsilon}}=0$
yield
\begin{equation}
\partial_zw+\tau(\dddd\Phi)w=-\partial_z\left(\frac {g}{\tau\partial_z
\Phi}\right) + (\partial_{z}\chi_1)\widetilde
R_{\Phi,\tau}g\quad\mbox{in}\,\,\Omega,\quad
w\vert_{\partial\Omega}=0.
\end{equation}
Note that by (\ref{mika}) and the fact that $g\vert_{\mathcal H}=0$,
we have:
\begin{equation}\label{ona}
\left\vert \partial_z\left(\frac {g}{\partial_z \Phi}\right)
\right\vert
= \left\vert \frac{\partial_zg}{\partial_z\Phi}
- \frac{g}{\partial_z\Phi}\frac{\partial_z^2\Phi}{\partial_z\Phi}
\right\vert \le \frac{C}{\Pi_{k=1}^\ell\vert x-\widetilde x_k\vert}.
\end{equation}
Consider the cut off function $\chi\in C_0^\infty(\Omega)$ such that
$$
\chi\ge 0, \quad\chi\vert_{B(0,\frac 12)}=1. $$ By (\ref{ona}) and
Proposition \ref{Proposition 3.0} B),
\begin{equation}\label{mir}
\widetilde R_{\Phi,\tau}\left(\sum_{k=1}^\ell
\chi((x-\widetilde x_k)\ln \vert\tau\vert)
\partial_z\left(\frac
{g}{\partial_z \Phi}\right)\right)\rightarrow 0
\quad \mbox{in}\quad
L^2(\Omega)\,\,\mbox{as}\,\, \vert\tau\vert\rightarrow +\infty.
\end{equation}
In fact, fixing large $\vert \tau\vert$, small $\delta>0$ and
$p >1$ such that $p-1$ is sufficiently small, we apply
Proposition \ref{Proposition 3.0} B) and (\ref{ona}) to conclude
\begin{eqnarray*}
&&\left\Vert \widetilde{R}_{\Phi,\tau}\left(
\sum_{k=1}^{\ell} \chi((x-\widetilde{x}_k)\ln \vert \tau\vert)
\partial_z \left( \frac{g}{\partial_z\Phi}\right)\right)
\right\Vert_{L^2(\Omega)}^2\\
&\le& C\sum_{k=1}^{\ell}\int_{B(\widetilde{x}_k,\delta)}
\vert \chi((x-\widetilde{x}_k)\ln \vert \tau\vert)\vert^p
\left\vert \partial_z \left( \frac{g}{\partial_z\Phi}\right)
\right\vert^p dx\\
&\le& C'\sum_{k=1}^{\ell}\int_{B(\widetilde{x}_k,\delta)}
\vert \chi((x-\widetilde{x}_k)\ln \vert \tau\vert)\vert^p
\frac{1}{\vert x-\widetilde{x}_k\vert^p} dx
\le C''\int^{\delta}_0 \vert \chi(\rho\ln\vert\tau\vert)\vert^p
\rho^{1-p}d\rho.
\end{eqnarray*}
Thus we get (\ref{mir}) by Lebesgue's theorem.

By Proposition \ref{Proposition 3.3}, we obtain
\begin{equation}\label{mir1}
\widetilde R_{\Phi,\tau}\left(\left(1-\sum_{k=1}^\ell
\chi((x-\widetilde x_k)\ln \vert\tau\vert)\right)
\partial_z\left(\frac {g}{\partial_z \Phi}\right)\right)\rightarrow
0\quad \mbox{in}\quad L^2(\Omega)\,\,\mbox{as}\,\,
\vert\tau\vert\rightarrow +\infty.
\end{equation}
Therefore (\ref{mir}) and (\ref{mir1}) yield
\begin{equation}\label{minaA}
\left\Vert \widetilde R_{\Phi,\tau}\left(\partial_z\left(\frac
{g}{\partial_z \Phi}\right)\right)\right\Vert_{L^2(\Omega)}=o(1)\quad
\mbox{as}\,\,\tau\rightarrow +\infty.
\end{equation}
Denote $\widetilde w = w + \frac 1\tau\chi_1\widetilde
R_{\Phi,\tau}(\partial_z(\frac {g}{\partial_z \Phi})).$

By (\ref{minaA}), it suffices to prove
\begin{equation}\label{(3.4AAA)}
\Vert\widetilde w\Vert_{L^2(\Omega)} = o\left(\frac 1\tau\right)
\quad \mbox{as}\,\, \vert\tau\vert\rightarrow +\infty.
\end{equation}
In terms of (2.19) and (\ref{(3.2)}), observe that
\begin{equation}\label{zanoza}
\partial_z\widetilde w + \tau(\partial_z \Phi)\widetilde w
= f\quad\mbox{in}\,\,\Omega, \quad \widetilde
w\vert_{\partial\Omega}=0,
\end{equation}
where $f = \frac{1}{\tau}(\partial_z\chi_1)\widetilde
R_{\Phi,\tau}(\partial_z(\frac {g}{\partial_z
\Phi}))+(\partial_{z}\chi_1)\widetilde R_{\Phi,\tau}g.$ By
(\ref{minaA}) and (\ref{opa}),
\begin{equation}\label{npz}
\Vert f\Vert_{L^2(\Omega)}=o\left(\frac 1\tau\right)\quad \mbox{as}\,\,
\vert\tau\vert\rightarrow +\infty.
\end{equation}
Noting $\widetilde{w}\vert_{\partial\Omega} = 0$, applying
Proposition \ref{Proposition 2.1} to equation (\ref{zanoza})
and using (\ref{npz}), we obtain (\ref{(3.4AAA)}).
As for the first term in (2.18), we can argue similarly.
The proof of the proposition is completed.
\end{proof}

\section{\bf Complex geometrical optics solutions}

In this section, we construct complex geometrical optics solutions for
the Schr\"odinger equation $\Delta+q_1$ with $q_1$ satisfying the
conditions of Theorem \ref{main}. Consider
\begin{equation}\label{2.1I}
L_1u=\Delta u+q_1u=0 \quad \text{in}\,\, \Omega.
\end{equation}


We will construct solutions to (\ref{2.1I}) of the form
\begin{equation}\label{mozila}
u_1(x)=e^{\tau\Phi(z)}(a(z)+a_0(z)/\tau)+e^{\tau\overline{\Phi(z)}}
\overline{(a(z)+a_1(z)/\tau)} +
e^{\tau\varphi}u_{11}+e^{\tau\varphi}u_{12}, \quad
u_1\vert_{\Gamma_0}=0.
\end{equation}

The function $\Phi$ satisfies (\ref{zzz}), (\ref{mika}) and
\begin{equation}\label{zzz1}
\mbox{Im}\, \Phi\vert_{\Gamma_0}=0.
\end{equation}

The amplitude function $a(z)$ is not identically zero
on $\overline \Omega$ and has the following properties:
\begin{equation}\label{im}
a\in C^2(\overline\Omega),\quad \partial_{\overline z}a \equiv 0,
\,\,\mbox{Re}\, a\vert_{\Gamma_0}=0.
\end{equation}
The function $u_{11}$ is given by
\begin{equation}
u_{11}=-\frac 14 e^{i\tau\psi}\widetilde
R_{\Phi,\tau}(e_1(\partial^{-1}_{\overline z} (aq_1)-M_1(z)))
- \frac 14 e^{-i\tau\psi}R_{\Phi,-\tau}(e_1(\partial^{-1}_{z}
(\overline{a(z)}q_1)-M_3(\overline{z})))
\end{equation}
\begin{eqnarray*}
&&-\frac {e^{i\tau\psi}}{\tau} \frac{e_2(\partial^{-1}_{\overline z}
(aq_1)-M_1(z))}{4\partial_z\Phi}
- \frac{e^{-i\tau\psi}}{\tau}\frac{e_2(\partial^{-1}_{z}
(\overline{a(z)}q_1)-M_3(\overline{z}))}{4\overline{\partial_z\Phi}}\\
&=& w_1e^{-\tau\varphi} + w_2e^{-\tau\varphi},
\end{eqnarray*}
where the polynomials $M_1(z)$ and $M_3(\bar z)$ satisfy
\begin{equation}\label{popsa}
\partial_z^j(\partial^{-1}_{\overline z}(aq_1)-M_1(z))=0,
\quad x\in \mathcal H, \thinspace j=0,1,2,
\end{equation}
\begin{equation}\label{popsa2}
\partial_{\overline z}^j(\partial^{-1}_{z}(\overline{a}q_1)(z)
- M_3(\overline z)) = 0, \quad x \in \mathcal H,\,\, j=0,1,2.
\end{equation}
The functions $e_1,e_2\in C^\infty(\Omega)$ are constructed so that
$e_1+e_2\equiv 1$ on $\overline{\Omega}$, $e_2$ vanishes in some
neighborhood of $\mathcal H$ and $e_1$ vanishes in a neighborhood of
$\partial\Omega$  and we set
$$
w_1= -\frac 14 e^{\tau\Phi}\widetilde
R_{\Phi,\tau}(e_1(\partial^{-1}_{\overline z} (aq_1)-M_1(z)))
- \frac 14 e^{\tau\overline{\Phi}}R_{\Phi,-\tau}(e_1(\partial^{-1}_{z}
(\overline{a(z)}q_1)-{M_3(\overline z)}))
$$
and
$$
w_2 = -\frac {e^{\tau\Phi}}{\tau}
\frac{e_2(\partial^{-1}_{\overline z}(aq_1)-M_1(z))}{4\partial_z\Phi}
- \frac{e^{\tau\overline{\Phi}}}{\tau}\frac{e_2(\partial^{-1}_{z}
(\overline{a(z)}q_1)-{M_3(\overline z)})}
{4\overline{\partial_z\Phi}}.
$$
Finally $a_0,a_1$ are holomorphic functions such that $$
(a_0(z)+\overline{a_1(z)})\vert_{\Gamma_0}=\frac{(\partial^{-1}_{\overline
z} (aq_1)-M_1(z))}{4\partial_z\Phi} +\frac{(\partial^{-1}_{z}
(\overline{a(z)}q_1)-M_3(\overline{z}))}{4\overline{\partial_z\Phi}}
$$
Then, noting $\partial_{\overline{z}}\overline{\Phi} =
\overline{\partial_z\Phi}$, (\ref{(3.2)}) and (\ref{(3.3)}), we have
\begin{eqnarray*}
&&\Delta w_1=4\partial_{z}\partial_{\overline z}w_1\\
&=& - \partial_{\overline z}(e^{\tau\Phi}\partial_{ z}\widetilde
R_{\Phi,\tau}(e_1(\partial^ {-1}_{\overline z} (aq_1)-M_1(z)))
+(\tau\partial_z\Phi)e^{\tau\Phi}\widetilde R_{\Phi,\tau}
(e_1(\partial^{-1}_{\overline z} (aq_1)-M_1(z)))\\
&-& \partial_{ z}(e^{\tau\overline\Phi}\partial_{\overline
z}R_{\Phi,-\tau}(e_1(\partial^{-1}_{ z} (\overline aq_1)-{M_3(\overline
z)}))+(\tau\overline{\partial_z\Phi})e^{\tau\overline\Phi}R_{\Phi,-\tau}
(e_1(\partial^{-1}_{z} (\overline aq_1)-{M_3(\overline z)}))\\
&=& -\partial_{\overline z}(e^{\tau\Phi}e_1(\partial^{-1}
_{\overline z} (aq_1)-M_1(z)))
- \partial_{z}(e^{\tau\overline\Phi}e_1(\partial^{-1}_{ z} (\overline
aq_1)-{M_3(\overline z)})).
\end{eqnarray*}
Moreover
\begin{eqnarray*}
&&\Delta w_2=4\partial_{z}\partial_{\overline z}w_2\\
&=& -\partial_{\overline
z}(e^{\tau\Phi}(e_2(\partial^{-1}_{\overline z}
(aq_1)-M_1(z)))-\partial_{z}(e^{\tau\overline\Phi}
e_2(\partial^{-1}_{ z} (\overline aq_1)-M_3(\overline z)))\\
&-& e^{\tau\Phi}\Delta \left(
\frac{e_2(\partial^{-1}_{\overline z}(aq_1)-M_1(z))}
{4\tau\partial_z\Phi}\right)
- e^{\tau\overline\Phi}\Delta \left(\frac{e_2(\partial^{-1}_{z}
(\overline{a(z)}q_1)-{M_3(\overline z)})}
{4\tau\overline{\partial_z\Phi}}\right).
\end{eqnarray*}
Therefore
\begin{eqnarray}\label{popa}
\Delta (u_{11}e^{\tau\varphi})
= \Delta (w_1+w_2)=-aq_1e^{\tau\Phi}
- \overline aq_1e^{\tau\overline\Phi}\\
- e^{\tau\Phi}\Delta\left(
\frac{e_2(\partial^{-1}_{\overline z}
(aq_1)-M_1(z))}{4\tau\partial_z\Phi}\right)
- e^{\tau\overline\Phi}\Delta\left(
\frac{e_2(\partial^{-1}_{z}(\overline{a(z)}q_1)
- M_3(\overline{z}))}{4\tau\overline{\partial_z\Phi}}\right).\nonumber
\end{eqnarray}
By (\ref{im}) and (\ref{zzz1}), observe that
\begin{equation}\label{ebtvoyu}
(e^{\tau\Phi(z)}a(z)+e^{\tau\overline{\Phi(z)}}\overline{a(z)})
\vert_{\Gamma_0}=0.
\end{equation}

By Proposition \ref{Proposition 2.3},  the
inhomogeneous problem
\begin{equation}\label{(3.18)}
\Delta (u_{12}e^{\tau\varphi})
+q_1u_{12}e^{\tau\varphi}=-q_1u_{11}e^{\tau \varphi}+h_1e^{\tau
\varphi}\quad \mbox{in}\,\,\Omega,
\end{equation}
\begin{equation}\label{(3.188)}
u_{12}=-u_{11} \quad \mbox{on $\Gamma_0$},
\end{equation}
has a solution
where
\begin{eqnarray}\label{moma}
h_1=e^{\tau i\psi}\Delta \left(\frac{e_2(\partial^{-1}_{z}
({a(z)}q_1)-{M_1(z)})}{4\tau{\partial_z\Phi}}\right)
+ e^{-\tau i\psi}\Delta\left(
\frac{e_2(\partial^{-1}_{z}(\overline{a(z)}q_1)-M_3(\overline{z}))}
{4\tau\overline{\partial_z\Phi}}\right)\nonumber\\
-a_0q_1e^{\tau\Phi}/\tau-\overline{a_1}q_1e^{\tau\overline{\Phi}}/\tau.
\end{eqnarray}
Then, by (\ref{im}) and (\ref{popa}) - (\ref{moma}), we see that
(\ref{2.1I}) is satisfied.

By Proposition \ref{Proposition 2.3} there exists a positive
$\tau_0$ such that for all $\vert \tau\vert>\tau_0$ there exists  a
solution to (\ref{(3.18)}), (\ref{(3.188)}) satisfying
\begin{equation}\label{(3.20)}
\Vert u_{12}\Vert_{L^2(\Omega)}= o(\frac{1}{\tau}).
\end{equation}
This can be done because
$$
\Vert q_1u_{11}+h_1\Vert_{L^2(\Omega)}\le
C(\delta)/\tau^{1-\delta}\quad \forall\delta\in(0,1);\,\,\Vert
u_{11}\Vert_{L^2(\partial\Omega)} = o(\frac{1}{\tau})
$$
and $(\nabla \varphi, \nu) = 0$ on $\Gamma_0.$ The latter fact
can be seen as
follows: On $\partial\Omega $, the Cauchy-Riemann equations imply
$$
(\nabla\varphi,\nu) = \nu_1\partial_{x_1}\varphi +
\nu_2\partial_{x_2}\varphi = \nu_1\partial_{x_2}\psi -
\nu_2\partial_{x_1}\psi=\frac{\partial\psi}{\partial\vec \tau},
$$
which is the tangential derivative of $\psi = \mbox{Im}\thinspace
\Phi$ on $\partial\Omega$.  By (\ref{zzz1}) we see that the tangential
derivative of $\psi$ vanishes on $\Gamma_0$.

Consider the Schr\"odinger equation
\begin{equation}\label{2.1II}
L_2v=\Delta v+q_2v=0 \quad \text{in}\,\, \Omega.
\end{equation}

We will construct solutions to (\ref{2.1II}) of the form
\begin{equation}\label{mozilaa}
v(x)=e^{-\tau{\Phi(z)}}{(a(z)+b_0(z)/\tau)}
+e^{-\tau\overline{\Phi(z)}} \overline{(a(z)+b_1(z)/\tau)}
+e^{-\tau\varphi}v_{11}+e^{-\tau\varphi}v_{12}, \quad v\vert_{
\Gamma_0}=0.
\end{equation}

The construction of $v$ repeats the corresponding steps of the
construction of $u_1.$ The only difference is that instead of $q_1$
and $\tau$, we use $q_2$ and $-\tau$ respectively.
We provide the details for
the sake of completeness. The function $v_{11}$ is given by
\begin{equation}
v_{11}=-\frac 14 e^{-i\tau\psi}\widetilde
R_{\Phi,-\tau}(e_1(\partial^{-1}_{\overline z} (q_2
a(z))-M_2(z)))
- \frac 14 e^{i\tau\psi}R_{\Phi,\tau}(e_1(\partial^{-1}_{z}
(q_2\overline{a(z)}) - M_4(\overline{z})))
\end{equation}
$$
+ \frac {e^{-i\tau\psi}}{\tau} \frac{e_2(\partial^{-1}_{\overline z}
(aq_2)-M_2(z))}{4\partial_z\Phi}
+ \frac{e^{i\tau\psi}}{\tau}\frac{e_2(\partial^{-1}_{z}
(\overline{a(z)}q_2)-M_4(\overline{z}))}{4\overline{\partial_z\Phi}},
$$
where
\begin{equation}\label{popsa1}
\partial_z^j(\partial^{-1}_{\overline z}(aq_2)-M_2(z)) = 0,
\quad x \in \mathcal H, \,\, j=0,1,2,
\end{equation}
\begin{equation}\label{popsa4}
\partial_{\overline z}^j(\partial^{-1}_{z}
(\overline{a}q_2)(z)-M_4(\overline{z})) = 0, \quad x \in \mathcal
H,\,\, j=0,1,2.
\end{equation}
Finally $b_0,b_1$ are holomorphic functions such that
$$(b_0+\overline b_1)\vert_{\Gamma_0}=-\frac{(\partial^{-1}_{\overline z}
(aq_2)-M_2(z))}{4\partial_z\Phi} - \frac{(\partial^{-1}_{z}
(\overline{a(z)}q_2)-M_4(\overline{z}))}{4\overline{\partial_z\Phi}}.
$$
Denote
\begin{eqnarray}
h_2=e^{-\tau i\psi}\Delta \left(\frac{e_2(\partial^{-1}_{z}
({a(z)}q_2)-{M_2(z)})}{4\tau{\partial_z\Phi}}\right)
+ e^{\tau i\psi}\Delta\left(
\frac{e_2(\partial^{-1}_{z} (\overline{a(z)}q_2)
- M_4(\overline z))}{4\tau\overline{\partial_z\Phi}}\right)\nonumber\\
-\frac{b_0(z)}{\tau}e^{-\tau\Phi(z)}-\frac{\overline{b_1(z)}}{\tau}e^{-\tau\overline{\Phi(z)}}.\nonumber
\end{eqnarray}

The function $v_{12}$ is a solution to the problem:
\begin{equation}\label{(3.18')}
\Delta (v_{12}e^{-\tau\varphi})
+q_2v_{12}e^{-\tau\varphi}=-q_2v_{11}e^{-\tau
\varphi}-h_2e^{-\tau\varphi}\quad \mbox{in}\,\,\Omega,
\end{equation}
\begin{equation}
v_{12}\vert_{\Gamma_0}=-v_{11}\vert_{\Gamma_0}.
\end{equation}
such that
\begin{equation}\label{(3.20')}
\Vert v_{12}\Vert_{L^2(\Omega)} = o(\frac{1}{\tau}).
\end{equation}
\section{\bf Proof of the theorem.}
 We first apply stationary phase with a general phase function $\Phi$ and then
we construct an appropriate weight function.

\begin{proposition}\label{Lemma 3.1}
Suppose that $\Phi$ satisfies (\ref{zzz}),(\ref{mika}) and
(\ref{zzz1}). Let $\{\widetilde x_1,\dots,\widetilde x_\ell\}$ be
the set of critical points of the function $Im\Phi$. Then for any
potentials $q_1,q_2\in C^{1+\alpha}(\overline{\Omega})$, $\alpha>0$ with the
same Dirichlet-to-Neumann maps and for any holomorphic function $a$
satisfying (\ref{im}), we have
\begin{equation}\label{nonsence}
2\sum_{k=1}^\ell\frac{\pi(q\vert a\vert^2) (\widetilde
x_k)\mbox{Re}\,e^{2i\tau\mbox{Im}\Phi(\widetilde{x_k})}}
{\vert(\mbox{det} \thinspace \mbox{Im}\Phi'')(\widetilde{x_k})
\vert^\frac 12}+\int_\Omega q(a_0b_0+\bar a_1\bar b_1)dx
\end{equation}
\begin{eqnarray*}
&+& \frac{1}{4}\int_{\Omega} \left( qa \frac{\partial_{\overline
z}^{-1}(aq_2)-M_2(z)} {\partial_z\Phi}+
q\overline{a}\frac{\partial_{z}^{-1}(q_2\overline{a})
-{M_4(\overline z)}}{\overline{\partial_z\Phi}}\right)dx\\
&-& \frac{1}{4}\int_\Omega\left( qa\frac{(\partial_{\overline
z}^{-1}(aq_1)-M_1(z))}{\partial_z\Phi} + q\overline
a\frac{(\partial_{z}^{-1}(\overline a q_1)-{ M_3(\overline
z)})}{\overline{\partial_z\Phi}}\right)dx=0,\quad \tau > 0,
\end{eqnarray*}
where we set
$$
q = q_1 - q_2.
$$
\end{proposition}
\begin{proof}
We note by the Cauchy-Riemann equations that $\{
\widetilde{x}_{1,1}+i\widetilde{x}_{1,2}, ...,
\widetilde{x}_{\ell,1}+i\widetilde{x}_{\ell,2} \} = \{ z \in
\overline{\Omega}\vert \thinspace \dddd\mbox{Im}\thinspace \Phi(z) =
0 \}$. Let $u_1$ be a solution to (\ref{2.1I}) and satisfy
(\ref{mozila}), and $u_2$ be a solution to the following equation
$$
\Delta u_2+q_2u_2=0\quad \mbox{in}\,\,\Omega,\quad
u_2\vert_{\partial \Omega}=u_1\vert_{\partial\Omega}.
$$
Since the Dirichlet-to-Neumann maps are equal, we have
$$
\nabla u_2 = \nabla u_1 \quad \mbox{on $\widetilde\Gamma$}.
$$
Denoting $u=u_1-u_2$, we obtain
\begin{equation}\label{(3.22)}
\Delta u+q_2u=-qu_1\quad \mbox{in}\,\,\Omega,\quad u\vert_{\partial
\Omega}=\frac{\partial u}{\partial \nu}\vert_{\widetilde\Gamma}=0.
\end{equation}

Let $v$ satisfy (\ref{2.1II}) and (\ref{mozilaa}).  We multiply
(\ref{(3.22)}) by $v$, integrate over $\Omega$ and we use
$v\vert_{\Gamma_0} = 0$ and $\frac{\partial u}{\partial\nu} = 0$ on
$\widetilde{\Gamma}$ to obtain $\int_{\Omega} qu_1v dx = 0$. By
(\ref{mozila}), (\ref{(3.20)}), (\ref{mozilaa}) and (\ref{(3.20')}),
we have
\begin{eqnarray}\label{(3.31)}
0 = \int_\Omega qu_1vdx=\int_\Omega q(a^2 + \overline a^2
+ \vert a\vert^2e^{\tau(\Phi-\overline \Phi)}
+ \vert a\vert^2e^{\tau(\overline \Phi-\Phi)} \nonumber\\
+\frac{a_0b_0}{\tau}+\frac{\overline{a_1b_1}}{\tau}+
u_{11}e^{\tau\varphi}(ae^{-\tau\Phi}+\overline
ae^{-\tau\overline\Phi})\nonumber\\+(ae^{\tau\Phi}+\overline
ae^{\tau\overline\Phi})v_{11}e^{-\tau\varphi})dx + o\left(\frac
1\tau\right), \quad \tau > 0.
\end{eqnarray}
The first and second terms in the asymptotic expansion of
(\ref{(3.31)}) are independent of $\tau$, so that
\begin{equation}
\int_\Omega q(a^2 + \overline a^2)dx=0.
\end{equation}

Using stationary phase (see p.215 in \cite{E}. cf. \cite{H}),
we obtain
\begin{equation}\label{za}
\int_\Omega q(\vert a\vert^2 e^{\tau(\Phi-\overline \Phi)} + \vert
a\vert^2e^{\tau(\overline \Phi-\Phi)})dx
= 2\sum_{k=1}^\ell\frac{\pi q\vert a\vert^2(\widetilde
x_k)\mbox{Re}\,e^{2\tau i\mbox{Im}\,\Phi(\widetilde x_k)}} {\tau
\vert(\mbox{det}\thinspace \mbox{Im}\Phi'')(\widetilde
x_k)\vert^\frac 12}+o\left(\frac{1}{\tau}\right).
\end{equation}
Here by the Cauchy-Riemann equations, we see that
sgn$(\mbox{Im}\thinspace \Phi''(\widetilde{x}_k)) = 0$, where
sgn $A$ denotes the signature of the matrix $A$, that is the number of positive eigenvalues of $A$ minus the number of negative
eigenvalues (e.g., \cite{E}, p.210).  Moreover we use (2.2) and the
Cauchy-Riemann equations to see that
$$
\mbox{det}\thinspace \mbox{Im} \thinspace \Phi''(z)
= -(\partial_{x_1}\partial_{x_2}\varphi)^2 - (\partial_{x_1}^2\varphi)^2
\ne 0
$$
since $\partial_z^2\Phi = -\frac{1}{2}\partial_{x_1}^2\varphi
- \frac{1}{2}i\partial_{x_1}\partial_{x_2}\varphi\ne 0$ in $\mathcal H$.
We calculate the two remaining terms
in (\ref{(3.31)}).  We get:
\begin{eqnarray}
&&\qquad \qquad \quad \int_\Omega qu_{11}e^{\tau\varphi}(ae^{-\tau\Phi}
+\overline ae^{-\tau\overline\Phi})dx \\
&=& -\frac{1}{4}\int_\Omega q\left\{ e^{\tau\Phi}\widetilde{
R}_{\Phi,\tau} (e_1(\partial_{\overline z}^{-1}(aq_1)-M_1(z))) \right.\nonumber\\
&+& \left.
e^{\tau\overline\Phi}R_{\Phi,-\tau}(e_1
(\partial_{z}^{-1}(\overline{a}q_1) - {M_3(\overline z)}))\right\}
(ae^{-\tau\Phi}+\overline ae^{-\tau\overline\Phi})dx \nonumber\\
&-& \int_\Omega\left(\frac {e^{\tau\Phi}}{\tau}
\frac{e_2(\partial^{-1}_{\overline
z}(aq_1)-M_1(z))}{4\partial_z\Phi} +
\frac{e^{\tau\overline\Phi}}{\tau}\frac{e_2(\partial^{-1}_{z}
(\overline{a(z)}q_1)-{M_3(\overline z)})}{4\overline{\partial_z\Phi}}
\right)q(ae^{-\tau\Phi}+\overline ae^{-\tau\overline\Phi})dx\nonumber\\
&=& -\frac{1}{4}\int_\Omega(qa\widetilde R_{\Phi,\tau}
(e_1(\partial_{\overline z}^{-1}(aq_1)-M_1(z)))
+ q\overline aR_{\Phi,-\tau}(e_1
(\partial_{z}^{-1}(\overline{a}q_1)-{M_3(\overline z)})))dx \nonumber\\
&-& \frac{1}{4}\int_\Omega(q\overline a\widetilde
R_{\Phi,\tau}(e_1(\partial_{\overline
z}^{-1}(aq_1)-M_1(z)))e^{\tau(\Phi-\overline\Phi)} +
qaR_{\Phi,-\tau}(e_1(\partial_{z}^{-1}(\overline aq_1)-{
M_3(\overline z)}))e^{-\tau(\Phi-\overline\Phi)})dx \nonumber\\
&-& \int_\Omega q\left(\frac {e^{\tau(\Phi-\overline\Phi)}}{\tau}
\frac{\overline ae_2(\partial^{-1}_{\overline z}(aq_1)-M_1(z))}
{4\partial_z\Phi} + \frac{e^{\tau(\overline\Phi-\Phi)}}{\tau}
\frac{a e_2(\partial^{-1}_{z}(\overline{a(z)}q_1) - {M_3(\overline
z)})}{4\overline{\partial_z\Phi}} \right)
dx\nonumber\\
&-& \int_\Omega q\left(\frac {a}{\tau}
\frac{e_2(\partial^{-1}_{\overline z}
(aq_1)-M_1(z))}{4\partial_z\Phi} + \frac {\overline
a}{\tau}\frac{e_2(\partial^{-1}_{z} (\overline{a(z)}q_1)-{M_3(\overline
z)})}{4\overline{\partial_z\Phi}}
\right)dx\nonumber\\
&\equiv& I_1+I_2+I_3+I_4.  \nonumber
\end{eqnarray}

We compute $I_1$ and $I_2$ separately. By Proposition \ref{Proposition
3.22}, (\ref{popsa}) and stationary phase (e.g., p.215 in
\cite{E}), we get
\begin{eqnarray}\label{i0}
&&\qquad \qquad \quad I_2 = -\frac{1}{4}\int_\Omega
(q\overline a\widetilde R_{\Phi,\tau}(e_1(\partial_{\overline z}^{-1}
(aq_1)-M_1(z))) e^{\tau(\Phi-\overline\Phi)} \\
&+& qaR_{\Phi,-\tau}(e_1(\partial_{z}^{-1}(\overline aq_1)-{
M_3(\overline z)}))e^{-\tau(\Phi-\overline\Phi)})dx\nonumber\\
&=& -\frac{1}{4}\int_\Omega \left(
\frac {e_1q\overline{a} }{\tau {\partial_z\Phi}}(\partial_{\overline
z}^{-1}(aq_1)-M_1(z)) e^{2i\tau\mbox{Im}\Phi} +\frac{
e_1qa}{\tau\overline{\partial_z\Phi}}
(\partial_z^{-1}(\overline{a}q_1)-{M_3(\overline z)})
e^{-2i\tau\mbox{Im}\Phi}\right) dx \nonumber\\
&+ &o\left(\frac{1}{\tau}\right) = o\left(\frac{1}{\tau}\right).
                                                       \nonumber
\end{eqnarray}

By Proposition \ref{Proposition 3.22}, we obtain

\begin{equation}\label{i1}
I_1= -\frac{1}{4\tau}\int_\Omega
e_1\left(qa\frac{(\partial_{\overline
z}^{-1}(aq_1)-M_1(z))}{\partial_z\Phi} + q\overline
a\frac{(\partial_{z}^{-1}(\overline aq_1) - {M_3(\overline
z)})}{\overline{\partial_z\Phi}}\right)dx + o\left(\frac
1\tau\right).
\end{equation}
Using stationary phase again and (\ref{popsa}) we conclude that
\begin{equation}\label{nona}
I_3 = o\left(\frac 1\tau\right).
\end{equation}
Similarly
\begin{eqnarray}
&&\qquad \quad \quad \qquad \quad
\int_\Omega qv_{11}e^{-\tau\varphi}(ae^{\tau\Phi}
+ \overline ae^{\tau\overline\Phi})dx\\
&=& -\frac{1}{4}\int_\Omega q\left\{ e^{-\tau\Phi}\widetilde
R_{\Phi,-\tau} (e_1(\partial_{\overline z}^{-1}(aq_2)-M_2(z))) \right.
\nonumber
\\&+&\left.
e^{-\tau\overline\Phi}R_{\Phi,\tau} (e_1(\partial_{z}^{-1}(\overline
aq_2)-{M_4(\overline z)}))\right\}
(ae^{\tau\Phi}+\overline ae^{\tau\overline\Phi})dx \nonumber\\
&+&\int_\Omega q\left(\frac {e^{-\tau\Phi}}{\tau}
\frac{e_2(\partial^{-1}_{\overline z}
(aq_2)-M_2(z)))}{4\partial_z\Phi} +
\frac{e^{-\tau\overline\Phi}}{\tau}\frac{e_2(\partial^{-1}_{z}
(\overline{a(z)}q_2)-{M_4(\overline z)})}{4\overline{\partial_z\Phi}}
\right)(ae^{\tau\Phi}+\overline ae^{\tau\overline\Phi})dx\nonumber\\
&=& -\frac{1}{4}\int_\Omega (qa\widetilde R_{\Phi,-\tau}(e_1
(\partial_{\overline z}^{-1}(aq_2)-M_2(z)))
+ q\overline aR_{\Phi,\tau}(e_1(\partial_{z}^{-1}(\overline
aq_2)-{ M_4(\overline z)})))dx \nonumber\\
&-& \frac{1}{4}\int_\Omega [q\overline ae^{\tau(\overline\Phi-\Phi)}
(\widetilde R_{\Phi,-\tau}(e_1(\partial_{\overline z}^{-1}
(aq_2)-M_2(z)))
+ qae^{\tau(\overline\Phi-\Phi)}R_{\Phi,\tau}
(e_1(\partial_{z}^{-1}(\overline aq_2)-{M_4(\overline z)}))]dx\nonumber\\
&+& \int_\Omega q\left(\frac {e^{-\tau(\Phi-\overline\Phi)}}{\tau}
\frac{\overline ae_2(\partial^{-1}_{\overline z}
(aq_2)-M_2(z)))}{4\partial_z\Phi} +
\frac{e^{\tau(\Phi-\overline\Phi)}}{\tau}\frac{ae_2(\partial^{-1}_{z}
(\overline{a(z)}q_2)-{M_4(\overline z)})}{4\overline{\partial_z\Phi}}
\right)dx                                             \nonumber\\
&+& \int_\Omega q\left(\frac {a}{\tau} \frac{e_2
(\partial^{-1}_{\overline z}(aq_2)-M_2(z)))}{4\partial_z\Phi} +
\frac {\overline a}{\tau}\frac{e_2(\partial^{-1}_{z}
(\overline{a(z)}q_2)-{M_4(\overline z)})}{4\overline{\partial_z\Phi}}
\right)dx                   \nonumber\\
&\equiv &J_1+J_2+J_3+J_4.\nonumber
\end{eqnarray}
By (\ref{popsa1}) and Proposition \ref{Proposition 3.22}, we have
\begin{equation}\label{j0}
J_1 = \frac{1}{4\tau}\int_\Omega e_1\left( qa\frac
{\partial_{\overline z}^{-1}(aq_2)-M_2(z)}{\partial_z\Phi} +
q\overline a\frac{\partial_{z}^{-1}(\overline{a}q_2) - {M_4(\overline
z)}}{\overline{\partial_z\Phi}}\right)dx + o\left(\frac
1\tau\right).
\end{equation}
The stationary phase argument, (\ref{popsa1}) and Proposition
\ref{Proposition 3.22} yield
\begin{equation}\label{j1}
J_2 = -\frac{1}{4}\int_\Omega [q\overline a
e^{\tau(\overline\Phi-\Phi)}\widetilde
R_{\Phi,-\tau}(e_1(\partial_{\overline z}^{-1}(aq_2)-M_2(z))) +
qae^{\tau(\overline\Phi-\Phi)}R_{\Phi,\tau}
(e_1(\partial_{z}^{-1}(\overline{a}q_2) - {M_4(\overline z)}))]dx =
o\left(\frac 1\tau\right).
\end{equation}
By the stationary phase argument and (\ref{popsa1}), we see that
\begin{equation}\label{nona1}
J_3 = o\left(\frac 1\tau\right).
\end{equation}
Therefore, applying (\ref{za}), (\ref{i0}), (\ref{j0}), (\ref{j1}),
(\ref{nona}) and (\ref{nona1}) in (\ref{(3.31)}), we conclude that
\begin{eqnarray}
&&2\sum_{k=1}^\ell\frac{\pi(q\vert a\vert^2)(\widetilde
x_k)\mbox{Re}\,e^{2i\tau\mbox{Im}\Phi(\widetilde
x_k)}}{\vert(\mbox{det} \thinspace \mbox{Im}\Phi'')(\widetilde
x_k)\vert^\frac 12}                                 \nonumber\\
&+& \frac{1}{4}\int_{\Omega} \left( qa \frac{\partial_{\overline
z}^{-1}(a(z)q_2)-M_2(z)} {\partial_z\Phi} +
q\overline{a}\frac{\partial_{z}^{-1}(q_2\overline{a(z)})
-{M_4(\overline z)}}{\overline{\partial_z\Phi}}\right)dx\nonumber\\
\qquad \qquad &&-\frac{1}{4}\int_\Omega\left(
qa\frac{\partial_{\overline z}^{-1}(q_1a)-M_1(z)}{\partial_z\Phi}
+q\overline a\frac{\partial_{z}^{-1}(q_1\overline a)-{ M_3(\overline
z)}}{\overline{\partial_z\Phi}}\right)dx=o(1).
\end{eqnarray}
as $\tau\rightarrow+\infty$.
Passing to the limit in this equality and applying Bohr's theorem
(e.g., \cite{BS}, p.393), we finish the proof of the proposition.
\end{proof}

%

 We need the following proposition in the construction of the phase function
$\Phi$.

Let
$\widetilde y_1,\dots, \widetilde y_m \in \Omega.$
Denote by $\mathcal R=(\mathcal R(\widetilde y_1),\dots, \mathcal
R(\widetilde y_m))$ the following operator:
$$
\mathcal R (\tilde y_k)g= (u(\widetilde
y_k),\partial_{z}u(\widetilde y_k),\partial^2_{z}u(\widetilde y_k)),
$$
where
\begin{equation}\label{Laplaces}
\partial_{\bar z}u=0\quad \mbox{in}\,\,\Omega,\quad \mbox{Im}\,u\vert_{\Gamma_0}=0,\quad
\mbox{Im}\,u\vert_{\widetilde\Gamma}=g.
\end{equation}

We have
\begin{proposition}\label{zopak} The operator
$\mathcal R:C^{\infty}_0(\widetilde{\Gamma}) \rightarrow {\Bbb C}^{3m}$ satisfies $\mbox{Im} \mathcal R={\Bbb C}^{3m}.$
\end{proposition}
\begin{proof}
We note that $\mbox{Im}\mathcal R  = {\Bbb
C}^{3m}$ if and only if the closure of $\mbox{Im}\mathcal R
$ is equal to ${\Bbb C}^{3m}.$ Our proof is by
contradiction.  Assume that
$$
\mbox{Im}\mathcal R \ne {\Bbb C}^{3m},
$$
then there exists a nonzero vector \begin{equation}\label{mazda}\vec
A=(A_0^1,A_1^1, A^1_2, \dots,A_0^m,A_1^m, A^m_2) \in {\Bbb
C}^{3m}\in  (\mbox{Im}\mathcal R
)^\perp.\end{equation} It is known
that the problem (\ref{Laplaces}) for a fixed $g$ has solution if and
only if
$$ \int_{\Gamma_0} g wdz=0$$ for all $w$ such that
\begin{equation} \label{ooo}\partial_{\bar z}w=0\quad \mbox{in}\,\Omega,\quad
w\vert_{\partial\Omega}= \overline{z'(s)}\gamma(s). \end{equation}
Here $\gamma(s)$ is a real-valued function, $z(s)$ is the
parametrization of $\partial\Omega.$

Let the function $p$ be a solution to the boundary value problem
\begin{equation}\label{mirroy}
\partial_{\bar z} p=\sum_{k=1}^m (A^k_0\delta(x-\widetilde
y_k)-A^k_1\partial_{ z}\delta(x-\widetilde y_k)+A^k_2\partial^2_{
z}\delta(x-\widetilde y_k)),
\end{equation}
\begin{equation}\label{mirroy1}\mbox{Re}[iz'(s)p]\vert_{\partial\Omega}=0.\end{equation} By (\ref{mazda}) a solution to
(\ref{mirroy}), (\ref{mirroy1}) exists. Let the function $u$ be a
solution to problem $(\ref{Laplaces})$. Since
$$
\int_{\partial\Omega} pg d z= (\mathcal Rg, \vec A)=0,
$$
 using the boundary condition (\ref{mirroy1}) and Holmgren's theorem, we have $p=\tilde p+w $ with $w$ solving
 (\ref{ooo}) and $\tilde p$ such that
$\supp \tilde p\subset \{\widetilde y_1,\dots, \widetilde y_m\}.$
Since $p$ is  a distribution we have that that
$p=\sum_{k=1}^m\sum_{\vert\alpha\vert \le j(k)} C_{k,\alpha}D^\alpha
\delta(x-\widetilde y_k).$
\\
This implies that
(\ref{mirroy}) is possible only if  $\vec A=0$  which is a
contradiction.
\end{proof}

\noindent{\bf End of proof of Theorem 1.1}

\begin{proof}
We will construct a complex geometrical optics solution of the form
(\ref{mozila}) where $\Phi$ and $a$ satisfy (\ref{zzz}), (\ref{mika}), (\ref{zzz1}) and
(\ref{im}).  Let $a(z)$ be a solution to the Riemann-Hilbert problem
$$
\partial_{\bar z} a=0\quad \mbox{in}\,\,\Omega, \quad \mbox{Re}\,
a\vert_{\Gamma_0}=0 $$ which is not identically zero in $\Omega.$
Let $\widehat x$ be an arbitrary point from $\Omega$ such that
$a(\widehat x) \ne 0.$

 Next we construct a holomorphic function $\Phi$ such that
$\widehat{x} \in \mathcal G \equiv \{ x\in \overline{\Omega}\vert
\thinspace
\partial_z \Phi(x) = 0\}$, $\mbox{Im}\thinspace
\Phi(\widehat x) \ne \mbox{Im}\thinspace \Phi(x)$ if $x \in \mathcal
G$ and $x \ne \widehat{x}$.



Now we construct the function   $\Phi.$ Let $\tilde \Omega$ be a bounded domain in $\R^2$ such that $\overline\Omega\subset \tilde\Omega$, $\Gamma_0\subset\partial\Omega$, $\partial\tilde\Omega\cap (\partial\Omega\setminus \Gamma_0)=\emptyset.$
 By Proposition \ref{zopak} there exists a holomorphic function $u$
in $\tilde \Omega$
such that $\mbox{Im}\,u\vert_{\Gamma_0}=0$ and $u(\widehat x)=\partial_{z}u
(\widehat x)=0,$ and $\partial_z^2 u(\widehat x) \ne 0$. In general
the function $u$ may have critical points on the boundary. Let $\Gamma_*\subset\partial\Omega\setminus\Gamma_0$  such that $u$ does not have any critical points on $\overline{\Gamma_*}.$

Next we construct  a holomorphic function $p$
such that $u+\epsilon p$ does not have a critical points on
$\partial\Omega$ for all sufficiently small positive $\epsilon$ and
$\mbox{Re}\,p\vert_{\Gamma_0}=0.$ In order to do this  we use Proposition \ref{MMM1} proven in the Appendix.  We set up appropriate Cauchy data for the Cauchy-Riemann equations (\ref{(4.112)}) (see Appendix). On $\Gamma_0$ we set $\mbox{Re}\,p\vert_{\Gamma_0}=0$ and $\frac{\partial p}{\partial \nu}=\frac{\partial \mbox{Im} p}{\partial\vec \tau}<0$ on $\overline{\Gamma_0}.$ Let $p$ be a holomorphic function in $\tilde \Omega$
such that $\mbox{Re}\,p\vert_{\Gamma_0}=0$ and $\frac{\partial p}{\partial \nu}=\frac{\partial \mbox{Im} p}{\partial\vec \tau}<0$ on $\overline{\Gamma_0}.$ Obviously the function $u+\epsilon p$ does not have any critical points on $\overline\Gamma_0$ for all nonzero $\epsilon.$ On the other hand it might have a critical points on
the remaining part of the boundary $\partial\Omega\setminus\overline{\Gamma_0}.$ The number of such a critical points is finite and the function $\vert\nabla u\vert^2$ has a zero of finite order at these points.  By using a conformal transformation if necessary, we may assume that  $\partial\Omega\setminus \Gamma_0$ is a segment on the line $\{x_2=0\}.$ Let $\{(y_k,0)\}_{k=1}^{\tilde N}$ be the set of critical points of the function $u$ on the boundary $\Gamma_0.$

We divide the set $\{y_k\}_{k=1}^{\tilde N}$ into two sets $\mathcal O_1$ and $\mathcal O_2$.
 Let us fix some point $y_k$. By Taylor's formula  $\frac{\partial u}{\partial x_1}(x_1,0)=c_1(x_1-y_k)^{\kappa_1+1}+o((x_1-y_k)^{\kappa_1+1})$ and $\frac{\partial u}{\partial x_2}(x_1,0)=c_2(x_1-y_k)^{\kappa_2+1}+o((x_1-y_k)^{\kappa_2+1})$ with some $(c_1,c_2)\ne 0.$
 If $c_2\ne 0$ and $ \kappa_2\le \kappa_1$  we say that $y_k\in \mathcal O_1$. If  $c_1\ne 0$ and $ \kappa_2> \kappa_1$ we say that $y_k\in \mathcal O_2$.

Let us consider two cases. Let $y_k\in \mathcal O_1.$ Then if $\kappa_2$ is odd we take the Cauchy data for a  holomorphic function $p_k$ be such that $Re p_k=0$ and $\frac{\partial Im p_k}{\partial \vec \tau}$ is positive near $y_k$ if $c_2$ is positive ,  $\frac{\partial Im p_k}{\partial\vec \tau}$ is negative  near $y_k$ if $c_2$ is negative   and small on $\partial \Omega\setminus \Gamma_*.$ If $\kappa_2$ is even and $\kappa_1\ne \kappa_2$ we take the Cauchy data such that  $\frac{\partial Im p_k}{\partial\vec \tau}(y_k)-1$, $\frac{\partial Re p_k}{\partial\vec \tau}(y_k)-1$ are small otherwise $\frac{1}{c_2}\frac{\partial Im p_k}{\partial\vec \tau}(y_k)\ne \frac{1}{c_1}\frac{\partial Re p_k}{\partial\vec \tau}(y_k).$

Let $y_k\in \mathcal O_2.$ Then if $\kappa_1$ is odd we take the holomorphic function $p_k$ such that  $\frac{\partial Rep_k}{\partial \vec \tau}$ is positive near $y_k$ if $c_1$ is positive ,  $\frac{\partial Re p_k}{\partial\vec \tau}$ is negative  near $y_k$ if $c_1$ is negative   and small on $\partial \Omega\setminus \Gamma_*.$ If $\kappa_1$ is even  we take $\frac{\partial Re p_k}{\partial\vec \tau}(y_k)-1,\frac{\partial Im p_k}{\partial\vec \tau}(y_k)-1$ to be small.
Now we have finished the  construction of a  Cauchy data on $\Gamma_0$ and in a neighborhood $ \mathcal U$ of the set $\{(y_k,0)\}_{k=1}^{\tilde N}.$  On the part
of the boundary $\partial\Omega\setminus (\Gamma_0\cup\mathcal U\cup \Gamma_*)$ we continue $\mbox{Im}\,p, \mbox{Re}\,p$ up to smooth functions.    By Proposition \ref{MMM1} and general results on a solvability of the boundary problem for $\partial_{\bar z}$ operator there exists a holomorphic function $p$   which satisfies the  above choice of the  Cauchy data.

Denote by $\mathcal H_\epsilon$
the set of critical points of the function $u+\epsilon p$ on
$\overline{\Omega}$. By the implicit function theorem, there exists
a neighborhood of $\widehat x$ such that for all small $\epsilon$ in
this neighborhood the function $u+\epsilon p$ has only one critical
point $\widehat x(\epsilon)$, this critical point is nondegenerate
and
\begin{equation} \label{mina}
\widehat x (\epsilon)\rightarrow \widehat x,\,\,\mbox{as}\,\, \quad \epsilon\rightarrow 0.
\end{equation}

Let us fix a sufficiently small $\epsilon.$ Let $\mathcal
H_{\epsilon} = \{ x_{k,\epsilon} \}_{1\le k \le N(\epsilon)}$. By
Proposition \ref{zopak}, there exists a holomorphic function $w$ such
that
\begin{equation}\label{zozo}
w\vert_{\Gamma_0}=0,\,\,w(x_{k,\epsilon})\ne w(x_{j,\epsilon})
\quad \mbox{for}\,\,\, k\ne j,\quad
\partial_z w\vert_{\mathcal H_{\epsilon}}=0,\quad
\partial^2_{z}w\vert_{\mathcal H_{\epsilon}}\ne 0.
\end{equation}
Denote $\Phi_\delta=u+\epsilon p+\delta w.$ For all sufficiently
small positive $\delta$
$$
\mathcal H_{\epsilon}\subset\mathcal
G_\delta=\{x\in\overline\Omega\vert \partial_z\Phi_\delta(x)=0\}
$$
and
\begin{equation}\label{ono}
\inf_{\forall y\in \mathcal
H_{\epsilon},\,\,\hat x(\epsilon)\ne y}\vert\Phi_\delta(\hat x(\epsilon))- \Phi_\delta(y) \vert>\hat C(\epsilon)>0, \quad C(\epsilon)=O(\delta).
\end{equation}
 We show now that for all small positive $\delta$ the critical
points of the function $\Phi_\delta$ are nondegenerate.
Let $\widetilde x$ be a  critical point of the function $u+\epsilon p.$
If $\widetilde x$ is a nondegenerate critical point, by the implicit
function theorem, there exists a ball $B(\widetilde x,\delta_1)$  such
that the function $\Phi_\delta$ in this ball has only one nondegenerate
critical point for all small $\delta.$ Let $\widetilde x$ be a
degenerate critical point of $u+\epsilon p.$
Without loss of generality we may assume that $\widetilde x=0$.
In some neighborhood of $0$, we have
$\partial_{z}\Phi_\delta=\sum_{k=1}^\infty c_k
z^{k+\hat k}-\delta \sum_{k=1}^\infty b_k z^k$ for some integer positive
$\hat k,$ some $c_1\ne 0$ and thanks to (\ref{zozo}) some $b_1\ne 0.$ Let $z_\delta\in \mathcal
G_\delta$ and $z_\delta\rightarrow 0.$ Then either \begin{equation}\label{noraa}
z_\delta=0\,\,\mbox{ or}\,\,\,
z_\delta^{\hat k}=\delta b_1/c_1+o(\delta).\end{equation}
Therefore $\partial^2_z\Phi(z_\delta)\ne 0$ for all small $\delta.$
Hence we can apply Proposition \ref{Lemma 3.1} to conclude
$$
\sum_{x\in \mathcal G_\delta}q(x)c(x)e^{2i\tau
\mbox{Im}\Phi_\delta(x)}=0.
$$
By (\ref{nonsence}) $c(\widehat x(\epsilon))$ is not equal to zero. Also we claim that for all small positive $\delta$
\begin{equation}\label{XXX}\mbox{Im}\Phi_\delta(\widehat x(\epsilon))\ne \mbox{Im}\Phi_\delta(x)\,\,\,
\forall x\in \mathcal G_\delta\,\,\mbox{ such that}\,\, \widehat x(\epsilon)\ne x.\end{equation}
 Really, suppose that there exists a sequence $\tilde x_\delta\in \mathcal G_\delta $ such that $\mbox{Im}\Phi_\delta(\widehat x(\epsilon))= \mbox{Im}\Phi_\delta(\tilde x_\delta)$ as $\delta\rightarrow +0.$ Then taking if it is necessary a subsequence  we have that $\tilde x_\delta\rightarrow \tilde x\in\mathcal H_\epsilon$ and $\hat x(\epsilon)\ne \tilde x.$ In that case $\mbox{Im}\Phi_\delta(\tilde x)=\mbox{Im}\Phi_\delta(\widehat x(\epsilon)).$
On the other hand since, $$\Phi_\delta=\Phi_\delta(\hat x(\epsilon))+\sum_{k=1}^\infty \frac{c_k}{k+\hat k+1}
(z-\tilde z)^{k+\hat k+1}-\delta \sum_{k=1}^\infty \frac {b_k}{k+1} (z-\tilde z)^{k+1}$$ we have from (\ref{noraa})
$$\Phi_\delta(\hat x(\epsilon))-\Phi_\delta(\hat x(\epsilon))=\Phi_\delta(\tilde x)-\Phi_\delta(\hat x(\epsilon))=\sum_{k=1}^\infty \frac{c_k}{k+\hat k+1}
(z_\delta-\tilde z)^{k+\hat k+1}-\delta \sum_{k=1}^\infty \frac {b_k}{k+1} (z_\delta-\tilde z)^{k+1}
$$
$$
=\sum_{k=1}^\infty \frac{c_k}{k+\hat k+1}
(z_\delta-\tilde z)^{k+1}(\delta b_1/c_1+o(\delta))-\delta \sum_{k=1}^\infty \frac {b_k}{k+1} (z_\delta-\tilde z)^{k+1}\ne 0
$$
for all sufficiently small positive $\delta.$

Since the exponents are linearly independent functions of $\tau$, thanks to (\ref{XXX})
we have $q(\widehat x(\epsilon))=0$.  Thus (\ref{mina}) implies
$q(\widehat x)=0$, finishing the proof.
\end{proof}

\section{\bf Appendix.}

Consider the Cauchy problem for the Cauchy-Riemann equations
\begin{equation}\label{(4.112)}
L(\phi,\psi)=(\frac{\partial \phi}{\partial x_1}-\frac{\partial
\psi}{\partial x_2},\frac{\partial \phi}{\partial
x_2}+\frac{\partial \psi}{\partial x_1}) =0\quad\mbox{in}\,\,\Omega,
\quad \left(\phi,\psi\right)\vert _{\Gamma_0} =B=(b_1(x),b_2(x)).
\end{equation}

The following proposition establishes the solvability of
(\ref{(4.112)}) for a dense set of Cauchy data.
\begin{proposition}\label{MMM1}
There exist a set $\mathcal O\subset C^1(\overline{\Gamma_0})$ such
that for each $B\in \mathcal O$ problem (\ref{(4.112)}) has at least
one solution $(\phi,\psi)\in C^2(\overline\Omega)$ and
$\overline{\mathcal O}=C^1(\overline{\Gamma_0}).$
\end{proposition}

\begin{proof}
Consider the following extremal problem
\begin{equation}\label{extr}
J(\phi,\psi)=\left\Vert (\phi,\psi) -B\right\Vert^2 _{H^2(\Gamma_0)}
+ \epsilon \Vert (\phi,\psi)\Vert_{H^2(\partial\Omega)}^2 + \frac
1\epsilon\left\Vert\Delta L(\phi,\psi)\right\Vert^2_{L^2(\Omega)}
\rightarrow \inf,
\end{equation}
\begin{equation}\label{extr1}
(\phi,\psi)\in \mathcal X.
\end{equation}
Here $\mathcal X=\left\{\delta(x)=(\delta_1,\delta_2)\vert \delta\in
H^2(\Omega), \Delta L \delta\in L^2(\Omega),
L\delta\vert_{\partial\Omega} =0,
\delta\vert_{\partial\Omega}\in H^2(\partial\Omega)\right\}$.

For each $\epsilon>0$ there exists a unique solution to
(\ref{extr}), (\ref{extr1}) which we denote as $(\widehat
\phi_\epsilon,\widehat \psi_\epsilon)$. By Fermat's theorem  (see e.g. \cite{AL} p. 155) we
have
$$
J'(\widehat \phi_\epsilon,\widehat \psi_\epsilon)[\delta]=0,\quad
\forall\delta\in\mathcal X.
$$

This equality can be written in the form
$$
\left((\widehat \phi_\epsilon,\widehat \psi_\epsilon)-B,
\delta\right) _{H^2(\Gamma_0)}+\epsilon ( (\widehat
\phi_\epsilon,\widehat
\psi_\epsilon),\delta)_{H^2(\partial\Omega)}+\frac 1\epsilon (\Delta
L(\widehat \phi_\epsilon,\widehat \psi_\epsilon),\Delta
L\delta)_{L^2(\Omega)}=0.
$$
This equality implies that the sequence $\{(\widehat
\phi_\epsilon,\widehat \psi_\epsilon)\}$ is bounded in
$H^2(\Gamma_0)$, the sequence $\{{\epsilon} (\widehat
\phi_\epsilon,\widehat \psi_\epsilon)\}$ converges to zero in
$H^2(\partial\Omega)$ and $\left\{\frac {1}{{\epsilon}} \Delta
L(\widehat \phi_\epsilon,\widehat \psi_\epsilon)\right\}$ is bounded
in $L^2(\Omega).$

Therefore there exist $q\in H^2(\Gamma_0)$ and $p\in L^2(\Omega)$
such that
\begin{equation}\label{zopa}
(\widehat \phi_\epsilon,\widehat \psi_\epsilon)-B \rightharpoonup q
\quad \mbox{ weakly in} \,\, H^2(\Gamma_0)
\end{equation}
and
\begin{equation}\label{zima}
\left(q,\delta \right) _{H^2(\Gamma_0)} + (p,\Delta
L\delta)_{L^2(\Omega)}=0\quad \forall \delta\in\mathcal X.
\end{equation}

Next we claim that
 \begin{equation}\label{Laplace}
 \Delta p=0\quad\mbox{in}\,\,\Omega
 \end{equation}
in the sense of distributions. Suppose that (\ref{Laplace}) is
already proved. This implies
 $$
 (p,\Delta L\delta)_{L^2(\Omega)}=0\quad \forall \delta\in H^4(\Omega), \quad L\delta\vert_{\partial\Omega}=
 \frac{\partial L\delta}{\partial\nu}\vert_{\partial\Omega}=0.
 $$
 This equality and (\ref{zima}) yield
 \begin{equation}\label{zima1}
\left(q,\delta \right) _{H^2(\Gamma_0)}=0\quad \forall \delta\in
H^4(\Omega), L\delta\vert_{\partial\Omega}=\frac{\partial
L\delta}{\partial\nu}\vert_{\partial\Omega}=0.
\end{equation}
Then using the trace theorem we conclude that $q=0$ and (\ref{zopa})
implies that
$$
(\widehat \phi_{\epsilon_k},\widehat
\psi_{\epsilon_k})-B\rightharpoonup 0\quad\mbox{weakly in}\,\,
{H^2(\Gamma_0)}.
$$
By the Sobolev embedding theorem
$$
(\widehat \phi_{\epsilon_k},\widehat \psi_{\epsilon_k})-B\rightarrow
0\quad\mbox{in}\,\, C^1(\overline{\Gamma_0}).$$ Therefore the
sequence $ \{(\widehat \phi_{\epsilon_k},\widehat
\psi_{\epsilon_k})-(\widetilde \phi_{\epsilon_k},\widetilde
\psi_{\epsilon_k})\}, $ with
$$L(\widetilde
\phi_{\epsilon_k},\widetilde \psi_{\epsilon_k})= L(\widehat
\phi_{\epsilon_k},\widehat
\psi_{\epsilon_k})\quad\mbox{in}\,\,\Omega,\,\,\,\widetilde
\psi_{\epsilon_k}\vert_{\Gamma_0}=0
$$
represents the desired approximation for the solution of the Cauchy
problem (\ref{(4.112)}).

Now we prove (\ref{Laplace}). Let $\widetilde x$ be an arbitrary
point in $\Omega$ and let $\widetilde \chi$ be a smooth function
such that it is zero in some neighborhood of
$\partial\Omega\setminus \Gamma_0$ and the set $\mathcal
B=\{x\in \Omega\vert \widetilde \chi(x)=1\}$ contains an open
connected subset $\mathcal F$ such that $\widetilde x\in \mathcal F$
and  $ \Gamma_0\cap \overline{\mathcal F}$ is an open set
in $\partial\Omega.$  By (\ref{zima})
$$
0=(p,\Delta L(\widetilde\chi\delta))_{L^2(\Omega)} =(\widetilde \chi
p,\Delta L \delta)_{L^2(\Omega)} +(p,[\Delta
L,\widetilde\chi]\delta)_{L^2(\Omega)}.
$$
That is,
\begin{equation}\label{vorona}
(\widetilde \chi p,\Delta L\delta)_{L^2(\Omega)}+([\Delta
L,\widetilde\chi]^*p,\delta)_{L^2(\Omega)}=0\quad \forall \delta\in
\mathcal X.
\end{equation}
This equality implies that $\widetilde\chi p\in H^1(\Omega).$

Next we take another smooth cut off function $\widetilde\chi_1$ such
that $\mbox{supp}\,\widetilde\chi_1\subset\mathcal B.$ A
neighborhood of $\widetilde x$ belongs to  $\mathcal B_1=\{x\vert
\widetilde\chi_1=1\}$, the interior of $\mathcal B_1$ is connected,
and $\mbox{ Int }\mathcal B_1\cap \mathcal P_\epsilon$ contains an
open subset $\mathcal O$ in $\partial\Omega.$  Similarly to
(\ref{vorona}) we have
\begin{equation}\label{vorona}
(\widetilde \chi_1 p,\Delta L\delta)_{L^2(\Omega)}+([\Delta
L,\widetilde\chi_1]^*p,\delta)_{L^2(\Omega)}=0.
\end{equation}

This equality implies that $\widetilde\chi_1 p\in H^2(\Omega).$ Let
$\omega$ be a domain such that $\omega\cap\Omega=\emptyset$,
$\partial\omega\cap\partial\Omega\subset\mathcal O$ contains an
open set in $\partial\Omega.$

We extend $p$ on $\omega$ by zero. Then
$$
(\Delta(\widetilde \chi_1
p),L\delta)_{L^2(\Omega\cup\omega)}+([\Delta
L,\widetilde\chi]^*p,\delta)_{L^2(\Omega\cup\omega)}=0.
$$
Hence
$$
L^*\Delta(\widetilde\chi_1 p)=0 \quad\mbox{in} \,\,\mbox{ Int
}\mathcal B_1\cup\omega,\quad p\vert_\omega=0.
$$
By  Holmgren's theorem $\Delta(\widetilde\chi_1 p)\vert _{\mbox
{Int }\mathcal B_1}=0$, that is, $(\Delta p)(\widetilde x)=0.$
\end{proof}

Now we prove a Carleman estimate for the Laplace operator.

\begin{proposition}\label{Theorem 2.1}Suppose that $\Phi$ satisfies (2.1), (2.2).
 Let  $u\in H^1_0(\Omega)$ is  a real valued function.  Then we have:

\begin{eqnarray}\label{suno4} \tau\Vert
ue^{\tau\varphi}\Vert^2_{L^2(\Omega)}+\Vert
ue^{\tau\varphi}\Vert^2_{H^1(\Omega)}+\Vert\frac{\partial
u}{\partial\nu}e^{\tau\varphi}\Vert^2_{L^2(\Gamma_0)}+
\tau^2\Vert\vert\frac{\partial\Phi}{\partial z} \vert
ue^{\tau\varphi}\Vert^2_{L^2(\Omega)}\nonumber
\\
\le C(\Vert (\Delta u)e^{\tau\varphi}\Vert^2_{L^2(\Omega)}+\tau
\int_{\tilde\Gamma}\vert \frac{\partial
u}{\partial\nu}\vert^2e^{2\tau\varphi}d\sigma).
\end{eqnarray}

\end{proposition}

\begin{proof}
Denote $\tilde v=ue^{\tau\varphi},\Delta u=f$ and $$ \Omega_+=\{x\in
\partial\Omega\vert (\nabla\varphi,\nu)>0\},\quad \Omega_-=\{x\in
\partial\Omega\vert (\nabla\varphi,\nu)<0\}.
$$
Observe that $\Delta=4\frac{\partial}{\partial
z}\frac{\partial}{\partial\bar z}$ and
$\varphi(x_1,x_2)=\frac{1}{2}(\Phi(z)+\overline{\Phi(z)})$.
Therefore
$$
e^{\tau\varphi}\Delta e^{-\tau\varphi} \tilde v=
(2\frac{\partial}{\partial z}-\tau\frac{\partial\Phi}{\partial
z})(2\frac{\partial}{\partial\bar z} -\tau\frac{\partial\bar\Phi}
{\partial \bar z})\tilde v=(2\frac{\partial}{\partial\bar z}
-\tau\frac{\partial\bar\Phi} {\partial \bar
z})(2\frac{\partial}{\partial z}-\tau\frac{\partial\Phi}{\partial
z})\tilde v=fe^{\tau\varphi}.
$$

Denote $\tilde w_1=\overline{Q(z)}(2\frac{\partial}{\partial\bar z}
-\tau\frac{\partial\bar\Phi} {\partial \bar z})\tilde v,\tilde
w_2=Q(z)(2\frac{\partial}{\partial z} -\tau\frac{\partial\Phi}
{\partial  z})\tilde v,\frac{\partial\Phi}{\partial
z}=\psi_1(x_1,x_2)+i\psi_2(x_1,x_2),$ $Q(z)$ is some holomorphic
function in $\Omega$ which does not have zeros in $\bar \Omega.$
Thanks to the zero Dirichlet boundary condition for $u$  we have $$
\tilde w_1\vert_{\partial\Omega}=2\overline{Q(z)}\partial_{\bar
z}\tilde
v\vert_{\partial\Omega}=(\nu_1+i\nu_2)\overline{Q(z)}\frac{\partial
\tilde v}{\partial \nu}\vert_{\partial\Omega},\,\, \tilde
w_2\vert_{\partial\Omega}=2Q(z)\partial_{ z}\tilde
v\vert_{\partial\Omega}=(\nu_1-i\nu_2)Q(z)\frac{\partial \tilde
v}{\partial \nu}\vert_{\partial\Omega}. $$

 By Proposition 2.6
\begin{eqnarray} \Vert ( \frac{\partial}{\partial
x_1}-i\psi_2\tau)\tilde w_1\Vert^2_{L^2(\Omega)}-\tau
\int_{\partial\Omega}(\nabla\varphi,\nu)\vert Q\vert^2\vert
\frac{\partial\tilde
v}{\partial\nu}\vert^2d\sigma+\mbox{Re}\int_{\partial\Omega}i((\nu_2
\frac{\partial}{\partial x_1}-\nu_1 \frac{\partial}{\partial
x_2})\tilde w_1)\overline{\tilde
w_1}d\sigma+\nonumber\\
+\Vert (i\frac{\partial}{\partial x_2}+\psi_1\tau)\tilde
w_1\Vert^2_{L^2(\Omega)}=\Vert
Qfe^{\tau\varphi}\Vert^2_{L^2(\Omega)}.\nonumber
\end{eqnarray}

and
\begin{eqnarray} \Vert (\frac{\partial}{\partial
x_1}+i\psi_2\tau)\tilde w_2\Vert^2_{L^2(\Omega)}-\tau
\int_{\partial\Omega}(\nabla\varphi,\nu)\vert Q\vert^2\vert
\frac{\partial\tilde
v}{\partial\nu}\vert^2d\sigma+\mbox{Re}\int_{\partial\Omega}i((-\nu_2
\frac{\partial}{\partial x_1}+\nu_1 \frac{\partial}{\partial
x_2})\tilde w_2)\overline{\tilde
w_2}d\sigma+\nonumber\\
+\Vert (i\frac{\partial}{\partial x_2}-\psi_1\tau)\tilde
w_2\Vert^2_{L^2(\Omega)}=\Vert
Qfe^{\tau\varphi}\Vert^2_{L^2(\Omega)}.\nonumber
\end{eqnarray}

Let us simplify the integral $\mbox{Re}i\int_{\partial\Omega}((\nu_2
\frac{\partial}{\partial x_1}-\nu_1 \frac{\partial}{\partial
x_2})\tilde w_1)\overline{\tilde w_1}d\sigma.$ We recall that
$\tilde v=ue^{\tau \varphi}$ and $\tilde
w_1=\overline{Q(z)}(\nu_1+i\nu_2)\frac{\partial \tilde v}{\partial
\nu}=\overline{Q(z)}(\nu_1+i\nu_2)\frac{\partial u}{\partial \nu}
e^{\tau \varphi}.$ Denote $A+iB=\overline{Q(z)}(\nu_1+i\nu_2).$ We get
\begin{eqnarray} \mbox{Re}\int_{\partial\Omega}i((\nu_2
\frac{\partial}{\partial x_1}-\nu_1 \frac{\partial}{\partial
x_2})\tilde w_1)\overline{\tilde
w_1}d\sigma=\nonumber\\
\mbox{Re}\int_{\partial\Omega}i((\nu_2 \frac{\partial}{\partial
x_1}-\nu_1 \frac{\partial}{\partial x_2})[(A+iB)\frac{\partial
u}{\partial\nu}
e^{\tau \varphi}])(A-iB)\frac{\partial u}{\partial\nu}e^{\tau\varphi}d\sigma=\nonumber\\
\mbox{Re}\int_{\partial\Omega}i[(\nu_2 \frac{\partial}{\partial
x_1}-\nu_1 \frac{\partial}{\partial x_2})(A+iB)] \vert\frac{\partial
\tilde v}{\partial\nu}\vert^2
)(A-iB)d\sigma+\nonumber\\
\mbox{Re}\int_{\partial\Omega}\frac i2(A^2+B^2)((\nu_2
\frac{\partial}{\partial x_1}-\nu_1 \frac{\partial}{\partial
x_2})\vert\frac{\partial \tilde v}{\partial\nu}\vert^2
d\sigma=\nonumber\\
\int_{\partial\Omega} (\partial_{\vec \tau} AB-\partial_{\vec\tau}
BA)\vert\frac{\partial \tilde v}{\partial\nu}\vert^2
d\sigma.\nonumber\end{eqnarray}

Now we simplify the integral
$\mbox{Re}\int_{\partial\Omega}i((-\nu_2 \frac{\partial}{\partial
x_1}+\nu_1 \frac{\partial}{\partial x_2})\tilde w_2)\overline{\tilde
w_2}d\sigma.$ We recall that $\tilde v=ue^{\tau \varphi}$ and
$\tilde w_2=(\nu_1-i\nu_2)Q(z)\frac{\partial \tilde v}{\partial
\nu}=2(\nu_1-i\nu_2)Q(z)\frac{\partial u}{\partial \nu} e^{\tau
\varphi}.$ A straightforward computation gives  \begin{eqnarray}
\mbox{Re}\int_{\partial\Omega}i((-\nu_2 \frac{\partial}{\partial
x_1}+\nu_1 \frac{\partial}{\partial x_2})\tilde w_2)\overline{\tilde
w_2}d\sigma=\nonumber\\
\mbox{Re}\int_{\partial\Omega}i((-\nu_2 \frac{\partial}{\partial
x_1}+\nu_1 \frac{\partial}{\partial x_2})[(A-iB)\frac{\partial
u}{\partial\nu}
e^{\tau \varphi}])(A+iB)\frac{\partial u}{\partial\nu}e^{\tau\varphi_1}d\sigma=\\
\mbox{Re}\int_{\partial\Omega}i[(-\nu_2 \frac{\partial}{\partial
x_1}+\nu_1 \frac{\partial}{\partial x_2})(A-iB)] \vert\frac{\partial
\tilde v}{\partial\nu}\vert^2
)(A+iB)d\sigma-\nonumber\\
\mbox{Re}\int_{\partial\Omega}\frac i2(A^2+B^2)((\nu_2
\frac{\partial}{\partial x_1}-\nu_1 \frac{\partial}{\partial
x_2})\vert\frac{\partial \tilde v}{\partial\nu}\vert^2
d\sigma=\nonumber\\
\int_{\partial\Omega}(\partial_{\vec\tau} AB-\partial_{\vec\tau} BA)
\vert\frac{\partial \tilde v}{\partial\nu}\vert^2
d\sigma.\nonumber\end{eqnarray} Using the above formula we obtain

\begin{eqnarray} \label{suno5}\Vert ( \frac{\partial}{\partial
x_1}+i\psi_2\tau)\tilde w_2\Vert^2_{L^2(\Omega)}+\Vert
(i\frac{\partial}{\partial x_2}-\psi_1\tau)\tilde
w_2\Vert^2_{L^2(\Omega)}-2\tau
\int_{\partial\Omega}(\nu,\nabla\varphi)\vert Q\vert^2\vert
\frac{\partial\tilde v}{\partial\nu}\vert^2d\sigma\nonumber\\
\Vert ( \frac{\partial}{\partial x_1}-i\psi_2\tau)\tilde
w_1\Vert^2_{L^2(\Omega)}+\Vert (i \frac{\partial}{\partial
x_2}+\psi_1\tau)\tilde w_1\Vert^2_{L^2(\Omega)}\nonumber
\\
+2\int_{\partial\Omega} (\partial_{\vec\tau} AB-\partial_{\vec\tau}
BA)\vert\frac{\partial \tilde v}{\partial\nu}\vert^2 d\sigma =
2\Vert Q fe^{\tau\varphi}\Vert^2_{L^2(\Omega)}.
\end{eqnarray}

Let $\tilde\psi_k$ be functions such that
$$
\frac{\partial\tilde \psi_1}{\partial
x_1}=\psi_2,\quad\frac{\partial\tilde \psi_2}{\partial
x_2}=\psi_1\quad \mbox{in}\,\,\Omega.
$$
We can rewrite equality (\ref{suno5}) in the form
\begin{eqnarray}\label{(2.20)} \Vert  \frac{\partial}{\partial
x_1}(e^{i\tilde\psi_1\tau}\tilde w_2)\Vert^2_{L^2(\Omega)}+\Vert
\frac{\partial}{\partial x_2}(e^{-i\tilde\psi_2\tau}\tilde
w_2)\Vert^2_{L^2(\Omega)}-2\tau
\int_{\partial\Omega}(\nu,\nabla\varphi)\vert Q\vert^2\vert
\frac{\partial\tilde v}{\partial\nu}\vert^2d\sigma\nonumber\\
\Vert  \frac{\partial}{\partial x_1}(e^{-i\tilde\psi_1\tau}\tilde
w_1)\Vert^2_{L^2(\Omega)}+\Vert \frac{\partial}{\partial
x_2}(e^{i\tilde\psi_2\tau}\tilde w_1)\Vert^2_{L^2(\Omega)}\nonumber
\\
+2\int_{\partial\Omega} (\partial_{\vec\tau} AB-\partial_{\vec\tau}
BA)\vert\frac{\partial \tilde v}{\partial\nu}\vert^2 d\sigma
={2}\Vert Q fe^{\tau\varphi}\Vert^2_{L^2(\Omega)}.
\end{eqnarray}

Observe that there exists some positive constant $C$ , independent of $\tau$,
such that \begin{eqnarray}\label{(2.21)}\frac
1C(\Vert\tilde w_1\Vert^2_{L^2(\Omega)}+\Vert\tilde
w_2\Vert^2_{L^2(\Omega)})\le \frac 12\Vert  \frac{\partial}{\partial
x_1}(e^{i\tilde\psi_2\tau}\tilde w_2)\Vert^2_{L^2(\Omega)}+\frac
12\Vert \frac{\partial}{\partial x_2}(e^{i\tilde\psi_1\tau}\tilde
w_2)\Vert^2_{L^2(\Omega)}\nonumber\\-\tau
\int_{\partial\Omega_-}(\nu,\nabla\varphi)\vert Q\vert^2\vert
\frac{\partial\tilde v}{\partial\nu}\vert^2d\sigma\nonumber\\
\frac 12\Vert  \frac{\partial}{\partial
x_1}(e^{-i\tilde\psi_1\tau}\tilde w_1)\Vert^2_{L^2(\Omega)}+\frac
12\Vert \frac{\partial}{\partial x_2}(e^{i\tilde\psi_2\tau}\tilde
w_1)\Vert^2_{L^2(\Omega)}.\end{eqnarray}

Since $\tilde v$ is a real-valued function we have
$$
\Vert 2\frac{\partial \tilde v}{\partial x_1}+\tau\psi_1\tilde
v\Vert^2_{L^2(\Omega)} +\Vert 2\frac{\partial \tilde v}{\partial
x_2}-\tau\psi_2\tilde v\Vert^2_{L^2(\Omega)}\le C_0(\Vert\tilde
w_1\Vert^2_{L^2(\Omega)}+\Vert\tilde w_2\Vert^2_{L^2(\Omega)}).
$$
Therefore
\begin{eqnarray}\label{(2.22)} 4\Vert \frac{\partial \tilde v}{\partial
x_1}\Vert^2_{L^2(\Omega)}
-2\tau\int_\Omega(\frac{\partial\psi_1}{\partial x_1}-\frac{\partial\psi_2}{\partial x_2})\tilde v^2dx\nonumber\\
+\Vert\tau\psi_1\tilde v\Vert^2_{L^2(\Omega)} +4\Vert \frac{\partial
\tilde v}{\partial x_2}\Vert^2_{L^2(\Omega)}+\Vert\tau\psi_2\tilde
v\Vert^2_{L^2(\Omega)} \le C_1(\Vert\tilde
w_1\Vert^2_{L^2(\Omega)}+\Vert\tilde w_2\Vert^2_{L^2(\Omega)}).
\end{eqnarray}

By the Cauchy-Riemann equations the second integral is zero.

Now since by assumption (\ref{mika}) the function $\Phi$ has zeros
of at most rder two we have
\begin{equation}\label{(2.23)}
\tau\Vert \tilde v\Vert^2_{L^2(\Omega)}\le C(\Vert \tilde
v\Vert^2_{H^1(\Omega)}+\tau^2\Vert\vert\frac{\partial\Phi}{\partial
z} \vert \tilde v\Vert^2_{L^2(\Omega)}).
\end{equation}

By (\ref{(2.22)}), (\ref{(2.23)})
\begin{equation}\label{(2.24)}
\tau\Vert \tilde v\Vert^2_{L^2(\Omega)}+\Vert \tilde
v\Vert^2_{H^1(\Omega)}+\tau^2\Vert\vert\frac{\partial\Phi}{\partial
z} \vert \tilde v\Vert^2_{L^2(\Omega)}\le C_1(\Vert\tilde
w_1\Vert^2_{L^2(\Omega)}+\Vert\tilde w_2\Vert^2_{L^2(\Omega)}).
\end{equation}
By (\ref{(2.24)}) we obtain from (\ref{(2.20)}), (\ref{(2.21)})
\begin{eqnarray} \frac{1}{C_5}(\tau\Vert \tilde
v\Vert^2_{L^2(\Omega)}+\Vert \tilde
v\Vert^2_{H^1(\Omega)}+\tau^2\Vert\vert\frac{\partial\Phi}{\partial
z} \vert\tilde v\Vert^2_{L^2(\Omega)})- \tau
\int_{\partial\Omega}(\nu,\nabla\varphi)\vert \frac{\partial\tilde
v}{\partial\nu}\vert^2d\sigma\nonumber
\\
+\int_{\partial\Omega}2(\partial_{\vec\tau} AB-\partial_{\vec\tau} BA)
\vert\frac{\partial \tilde v}{\partial\nu}\vert^2 d\sigma \le \Vert
fe^{\tau\varphi}\Vert^2_{L^2(\Omega)}+\tau \int_{\tilde
\Gamma}\vert(\nu,\nabla\varphi)\vert\vert \frac{\partial\tilde
v}{\partial\nu}\vert^2d\sigma.
\end{eqnarray}
Using Proposition \ref{MMM1} we make a choice of $Q(z)$ such
that $(\partial_{\vec\tau} AB-\partial_{\vec\tau} BA)$ is positive on
$\bar\Gamma_0.$

This  concludes the proof of the Proposition.
\end{proof}

%


\begin{thebibliography}{99} %

\bibitem{AL} V. Alekseev, V. Tikhomirov, S. Fomin, \textit{Optimal
Control}, Consultants Bureau, New York, 1987.

\bibitem{AP} K.~Astala, L.~P\"aiv\"arinta, \textit{Calder{\'o}n's
inverse conductivity problem in the plane}, Ann. of Math., \textbf{163}
(2006), 265--299.

\bibitem {ALP} K.\
Astala, M.\ Lassas, and L.\ P\"aiv\"airinta,
Calder\'on's inverse
problem for anisotropic conductivity in the plane,
{\it Comm.
Partial Diff. Eqns.} {\bf 30} (2005), 207--224.

\bibitem{BS} A.~B\"ottcher, B.~Silvermann, \textit{Analysis
of Toeplitz Operators},  Springer-Verlag, Berlin, 2006.

\bibitem{BT} R. Brown, R. Torres, \textit{Uniqueness in the inverse
conductivity problem for conductivities with $3/2$ derivatives in
$L^p, p>2n,$ } J. Fourier Analysis Appl {\bf 9} (2003), 1049--1056.

\bibitem{B-U} R. Brown, G. Uhlmann,  \textit{Uniqueness in the
 inverse conductivity problem with less regular conductivities in two
dimensions,} Comm. Partial Differential Equations, {\bf 22} (1997),
1009--1027.

\bibitem{Bu} A. Bukhgeim, \textit{Recovering the potential from Cauchy
data in two dimensions},  J. Inverse Ill-Posed Probl., {\bf 16} (2008),
19--34.

\bibitem{BuU}  A. Bukhgeim, G.\ Uhlmann, \textit{Recovering a potential
from partial Cauchy data}, Comm. Partial Differential Equations,
{\bf 27} (2002), 653--668.

\bibitem{C}  A. P.\ Calder\'on,   \textit{On an inverse boundary value
problem,} in \emph{Seminar on Numerical Analysis and its
Applications to Continuum Physics}, 65--73, Soc. Brasil. Mat., R\'io
de Janeiro, 1980.

\bibitem{ChengYama} J.\ Cheng, M.\ Yamamoto,
\textit{Determination of two convection coefficients from Dirichlet
to Neumann map in the two-dimensional case}, SIAM J. Math. Anal.,
{\bf35} (2004), 1371--1393.

\bibitem{DKSjU} D.~Dos Santos Ferreira, C.~Kenig, J.~Sj\"ostrand,
G.~Uhlmann, \textit{Determining a magnetic Schr\"odinger operator
from partial Cauchy data}, Comm. Math. Phys., \textbf{271} (2007),
467--488.

\bibitem{DKSaU} D. Dos Santos Ferreira, C. Kenig, M. Salo,
G. Uhlmann, \textit{Limiting Carleman weights and anisotropic
inverse problems}, preprint:  arXiv:0803.3508 .


\bibitem{E} L.~Evans, \textit{Partial Differential Equations},
American Mathematical Society, Providence, RI, 2000.

\bibitem {GLU}
A.\ Greenleaf, M.\ Lassas, and G.\ Uhlmann, \textit{The Calder\'on
problem for conormal potentials, I: Global uniqueness and
reconstruction},  Comm. Pure Appl. Math, {\bf 56} (2003), 328--352.

\bibitem{HW} H. Heck and J.-N. Wang, \textit{Stability estimates for
the inverse boundary value problem by partial
Cauchy data}, Inverse Problems, {\bf 22} (2006), 1787--1796.

\bibitem{H} L.~H\"ormander, \textit{The Analysis of Linear Partial
Differential Operators I}, Springer-Verlag, Berlin, 1985.

\bibitem{IGM} O.\ Imanuvilov, G. \ Uhlmann, M\ Yamamoto,
\textit {Partial data for the Calder\'on problem
in two dimensions}, preprint:  arXiv:0809.3037  (2008).

\bibitem{I} V. Isakov, \textit{On uniqueness in the inverse conductivity problem with local data}, Journal of Inverse Problems and Imaging, {\bf 1} (2008), 95-105.

\bibitem{KU} H. Kang, G. Uhlmann, \textit{Inverse problems for the Pauli Hamiltonian in two
dimensions}, Journal of Fourier Analysis and Applications, {\bf 10} (2004), 201-215.

\bibitem{KSU} C.~Kenig, J.~Sj\"ostrand, G.~Uhlmann,
\textit{The Calder\'on problem with partial data}, Ann. of Math.,
\textbf{165} (2007), 567--591.

\bibitem{K} K. Knudsen,\textit{ The Calder\'on problem with
partial data for less smooth conductivities}, Comm. Partial
Differential Equations, {\bf 31} (2006), 57--71.

\bibitem{KS} K.~Knudsen, M.~Salo,
\textit{Determining nonsmooth first order terms from partial
boundary measurements}, Inverse Problems and Imaging, \textbf{1}
(2007), 349--369.

\bibitem{K-V}R. Kohn  and M. Vogelius,  Identification of an
unknown conductivity by means of measurements at the boundary,
in Inverse Problems, edited by D. McLaughlin,
{\em SIAM-AMS Proceedings}, {\bf 14}(1984), 113-123.

\bibitem{N} A.~Nachman, \textit{Global uniqueness for a two-dimensional
inverse boundary value problem}, Ann. of Math., \textbf{143} (1996),
71--96.



\bibitem{PPU} L. P\"aiv\"arinta, A. Panchenko and G. Uhlmann,
\textit{Complex geometrical optics for Lipschitz conductivities,}
Revista Matematica Iberoamericana, {\bf 19} (2003), 57-72.


\bibitem{SuU} Z.\ Sun and G.\ Uhlmann, Anisotropic inverse
problems in two dimensions,
{\it Inverse Problems}, {\bf 19} (2003), 1001-1010.

\bibitem{S}J.\ Sylvester, An anisotropic inverse boundary value problem,
{\em Comm. Pure Appl. Math.},  {\bf 43}(1990), 201--232.

\bibitem{SU} J.~Sylvester, G.~Uhlmann, \textit{A global uniqueness
theorem for an inverse boundary value problem}, Ann. of Math.,
\textbf{125}  (1987), 153--169.


\bibitem{T} L. Tzou, Stability estimates for coefficients of magnetic
Schr\"odinger equation from full and partial measurements, to appear
Comm. Partial Differential Equations.

\bibitem{U} G.~Uhlmann, \textit{Commentary on Calder\'on's paper (29)
``On an inverse boundary value problem",} Selected Papers of A.P.
Calder\'on, edited by Alexandra Bellow, Carlos Kenig and Paul
Malliavin, AMS, (2008), 623--636.


\bibitem{VE}  I.~Vekua, \textit{Generalized Analytic Functions},
Pergamon Press, Oxford, 1962.
\end{thebibliography}
\end{document}